\documentstyle[11pt]{article}
  \newlength{\standardunitlength}
\newtheorem{cor}{Corollary} \newtheorem{lemma}{Lemma}
\newtheorem{theorem}{Theorem} 
\newenvironment{proof}{\noindent {\sc Proof:}}{$\Box$ \vspace{2 ex}}
\begin{document}

\begin{center} {\bf Random Matrix Theory over Finite Fields}
\end{center}

\begin{center}
By Jason Fulman
\end{center}

\begin{center}
Stanford University
\end{center}

\begin{center}
Department of Mathematics
\end{center}

\begin{center}
Building 380, MC 2125
\end{center}

\begin{center}
Stanford, CA 94305, USA
\end{center}

\begin{center}
http://math.stanford.edu/$\sim$ fulman
\end{center}

\begin{center}
fulman@math.stanford.edu
\end{center}

\begin{center}
First version: March 28, 2000
\end{center}

\begin{center}
Current version: March 23, 2001
\end{center}

\newpage
\tableofcontents

\section{Introduction} \label{Introduction}

	A natural problem is to understand what a typical element of
the finite general linear group $GL(n,q)$ ``looks like''. Many of the
interesting properties of a random matrix depend only on its conjugacy
class. The following list of questions one could ask are of this type:

\begin{enumerate}

\item How many Jordan blocks are there in the rational canonical form
of a random matrix?

\item What is the distribution of the order of a random matrix?

\item What is the probability that the characteristic polynomial of a
random matrix has no repeated factors?

\item What is the probability that the characteristic polynomial of a
random matrix is equal to its minimal polynomial?

\item What is the probability that a random matrix is semisimple
(i.e. diagonalizable over the algebraic closure $\bar{F_q}$ of the
field of $q$ elements)?

\end{enumerate} As Section \ref{cycleindex} will indicate, answers to
these questions have applications to the study of random number
generators, to the analysis of algorithms in computational group
theory, and to other parts of group theory. Section \ref{cycleindex}
describes a unified approach to answering such probability questions
using cycle index generating functions. As an example of its power, it
is proved independently in \cite{F1} and \cite{W2} that the $n
\rightarrow \infty$ limit of answer to Question 3 is
$(1-\frac{1}{q^5})/(1+\frac{1}{q^3})$. There is (at present) no other
method for deriving this result and generating functions give
effective bounds on the convergence rate to the limit. Extensions of
the cycle index method to the set of all matrices and to other finite
classical groups are sketched.

	Section \ref{GLSec} gives a purely probabilistic picture of
what the conjugacy class of a random element of $GL(n,q)$ looks
like. The main object of study is a probability measure $M_{GL,u,q}$
on the set of all partitions of all natural numbers. This measure is
connected with the Hall-Littlewood symmetric functions. Exploiting
this connection leads to several methods for growing random partitions
distributed as $M_{GL,u,q}$ and gives insightful probabilistic proofs
of group theoretic results. We hope to convince the reader that the
interplay between probability and symmetric functions is beautiful and
useful. A method is given for sampling from $M_{GL,u,q}$ conditioned
to live on partitions of a fixed size (which amounts to studying
Jordan form of unipotent elements) and for sampling from a $q$-analog
of Plancherel measure (which is related to the longest increasing
subsequence problem of random permutations).

	Section \ref{GLSec} goes on to describe a probabilistic
approach to $M_{GL,u,q}$ using Markov chains. This connection is
quite surprising, and we indicate how it leads to a simple and
motivated proof of the Rogers-Ramanujan identities. The measure
$M_{GL,u,q}$ has analogs for the finite unitary, symplectic, and
orthogonal groups. As this is somewhat technical these results are
omitted and pointers to the literature are given. However we remark
now that while the analogs of the symmetric function theory viewpoint
are unclear for the finite symplectic and orthogonal groups, the
connections with Markov chains carry over. Thus there is a coherent
probabilistic picture of the conjugacy classes of the finite classical
groups.

	Section \ref{triangular} surveys probabilistic aspects of
conjugacy classes in $T(n,q)$, the group of $n \times n$ upper
triangular matrices over the field $F_q$ with $1$'s along the main
diagonal. Actually a simpler object is studied, namely the Jordan form
of randomly chosen elements of $T(n,q)$. From work of Borodin and
Kirillov, one can sample from the corresponding measures on
partitions. We link their results with symmetric function theory and
potential theory on Bratteli diagrams.

	The field surveyed in this article is young and evolving. The
applications to computational group theory call for extensions of
probability estimates discussed in Section \ref{cycleindex} to maximal
subgroups of finite classical groups. It would be marvellous if the
program surveyed here carries over; this happens for the finite affine
groups \cite{F9}. The first step is understanding conjugacy classes
and partial results can be found in the thesis \cite{M}.
 	
	We close with a final motivation for the study of conjugacy
classes of random matrices over finite fields. The past few years have
seen an explosion of interest in eigenvalues of random matrices from
compact Lie groups. For the unitary group $U(n,C)$ over the complex
numbers, two matrices are in the same conjugacy class if and only if
they have the same set of eigenvalues. Hence, at least in this case,
which is related to the zeroes of the Riemann zeta function
\cite{KeSn}, the study of eigenvalues is the same as the study of
conjugacy classes.
		
	As complements to this article, the reader may enjoy the
surveys \cite{Sh2},\cite{Py1},\cite{Py2}. These articles discuss
probabilistic and enumerative questions in group theory and have
essentially no overlap with the program surveyed here.

\section{Cycle Index Techniques} \label{cycleindex}

	Before describing cycle index techniques for the finite
classical groups, we mention that the cycle index techniques here are
modelled on similar techniques for the study of conjugacy class
functions on the symmetric groups. For a permutation $\pi$, let
$n_i(\pi)$ be the number of length $i$ cycles of $\pi$. The cycle
index of a subgroup $G$ of $S_n$ is defined as \[ \frac{1}{|G|}
\sum_{\pi \in G} \prod_{i \geq 1} x_i^{n_i(\pi)} \] and is called a
cycle index because it stores information about the cycle structure of
elements of $G$. Applications of the cycle index to graph theory and
chemical compounds are exposited in \cite{PR}. It is standard to refer
to the generating function \[ 1+\sum_{n \geq 1} \frac{u^n}{n!}
\sum_{\pi \in S_n} \prod_{i \geq 1} x_i^{n_i(\pi)} \] as the cycle
index or cycle index generating function of the symmetric groups. From
the fact that there are $\frac{n!}{\prod_i n_i!  i^{n_i}}$ elements in
$S_n$ with $n_i$ cycles of length $i$, one deduces Polya's result that
this generating function is equal to $\prod_{m \geq 1}
e^{\frac{x_mu^m}{m}}$. This allows one to study conjugacy class
functions of random permutations (e.g. number of fixed points, number
of cycles, the order of a permutation, length of the longest cycle) by
generating functions. We refer the reader to \cite{Ko} for results in
this direction using analysis and to \cite{SL} for results about cycle
structure proved by a probabilistic interpretation of the cycle index
generating function.

	Subsection \ref{subsecGL} reviews the conjugacy classes of
$GL(n,q)$ and then discusses cycles indices for $GL(n,q)$ and
$Mat(n,q)$, the set of all $n \times n$ matrices with entries in the
field of $q$ elements. Subsection \ref{applic} describes applications
of cycle index techniques. Subsection \ref{Classical} discusses
generalizations of cycle indices to the finite classical groups.

	It is necessary to recall some standard notation. Let
$\lambda$ be a partition of some non-negative integer $|\lambda|$ into
integer parts $\lambda_1 \geq \lambda_2 \geq \cdots \geq 0$. We will
also write $\lambda \vdash n$ if $\lambda$ is a partition of $n$. Let
$m_i(\lambda)$ be the number of parts of $\lambda$ of size $i$, and
let $\lambda'$ be the partition dual to $\lambda$ in the sense that
$\lambda_i' = m_i(\lambda) + m_{i+1}(\lambda) + \cdots$. Let
$n(\lambda)$ be the quantity $\sum_{i \geq 1} (i-1) \lambda_i$ and let
$(\frac{u}{q})_i$ denote $(1-\frac{u}{q}) \cdots (1-\frac{u}{q^i})$.

\subsection{The General Linear Groups} \label{subsecGL}

	To begin we follow Kung \cite{Kun} in defining a cycle index
for $GL(n,q)$. First it is necessary to understand the conjugacy
classes of $GL(n,q)$. As is explained in Chapter 6 of the textbook
$\cite{Her}$, an element $\alpha \in GL(n,q)$ has its conjugacy class
determined by its rational canonical form. This form corresponds to
the following combinatorial data. To each monic non-constant
irreducible polynomial $\phi$ over $F_q$, associate a partition
(perhaps the trivial partition) $\lambda_{\phi}$ of some non-negative
integer $|\lambda_{\phi}|$. Let $deg(\phi)$ denote the degree of
$\phi$. The only restrictions necessary for this data to represent a
conjugacy class are that $|\lambda_z| = 0$ and $\sum_{\phi}
|\lambda_{\phi}| deg(\phi) = n.$

	An explicit representative of this conjugacy class may be
given as follows. Define the companion matrix $C(\phi)$ of a
polynomial $\phi(z)=z^{deg(\phi)} + \alpha_{deg(\phi)-1}
z^{deg(\phi)-1} + \cdots + \alpha_1 z + \alpha_0$ to be:

\[ \left( \begin{array}{c c c c c}
		0 & 1 & 0 & \cdots & 0 \\
		0 & 0 & 1 & \cdots & 0 \\
		\cdots & \cdots & \cdots & \cdots & \cdots \\
		0 & 0 & 0 & \cdots & 1 \\
		-\alpha_0 & -\alpha_1 & \cdots & \cdots & -\alpha_{deg(\phi)-1}
	  \end{array} \right) \] Let $\phi_1,\cdots,\phi_k$ be the polynomials such that
$|\lambda_{\phi_i}|>0$. Denote the parts of $\lambda_{\phi_i}$ by
$\lambda_{\phi_i,1} \geq \lambda_{\phi_i,2} \geq \cdots $. Then a matrix
corresponding to the above conjugacy class data is

\[ \left( \begin{array}{c c c c}
		R_1 & 0 & 0 &0 \\
		0 & R_2 & 0 & 0\\
		\cdots & \cdots & \cdots & \cdots \\
		0 & 0 & 0 & R_k
	  \end{array} \right) \] where $R_i$ is the matrix

\[ \left( \begin{array}{c c c}
		C(\phi_i^{\lambda_{\phi_i,1}}) & 0 & 0  \\
		0 & C(\phi_i^{\lambda_{\phi_i,2}}) & 0 \\
		0 & 0 & \cdots 
	  \end{array} \right) \]

	For example, the identity matrix has $\lambda_{z-1}$ equal to
$(1^n)$ and all other $\lambda_{\phi}$ equal to the emptyset. An
elementary transvection with $a \neq 0$ in the $(1,2)$ position, ones
on the diagonal and zeros elsewhere has $\lambda_{z-1}$ equal to
$(2,1^{n-2})$ and all other $\lambda_{\phi}$ equal to the
emptyset. For a given matrix only finitely many $\lambda_{\phi}$ are
non-empty. Many algebraic properties of a matrix can be stated in
terms of the data parameterizing its conjugacy class. For instance the
characteristic polynomial of $\alpha \in GL(n,q)$ is equal to
$\prod_{\phi} \phi^{|\lambda_{\phi}(\alpha)|}$ and the minimal
polynomial of $\alpha$ is equal to $\prod_{\phi}
\phi^{|\lambda_{\phi,1}(\alpha)|}$. Furthermore $\alpha$ is semisimple
(diagonalizable over the algebraic closure $\bar{F_q}$) precisely when
all $\lambda_{\phi}(\alpha)$ have largest part at most 1.

	To define the cycle index for $Z_{GL(n,q)}$, let
$x_{\phi,\lambda}$ be variables corresponding to pairs of polynomials
and partitions. Define

	\[ Z_{GL(n,q)} = \frac{1}{|GL(n,q)|} \sum_{\alpha \in GL(n,q)}
\prod_{\phi: |\lambda_{\phi}(\alpha)|>0} x_{\phi,\lambda_{\phi}(\alpha)}. \]
Note that the coefficient of a monomial is the probability of
belonging to the corresponding conjugacy class, and is therefore equal
to one over the order of the centralizer of a representative. It is
well known (e.g. easily deduced from page 181 of \cite{Mac}) that one
over the order of the centralizer of conjugacy class of $GL(n,q)$
corresponding to the data $\{\lambda_{\phi}\}$ is \[
\frac{1}{\prod_{\phi} q^{deg(\phi) \cdot \sum_i (\lambda_{\phi,i}')^2}
\prod_{i \geq 1} (\frac{1}{q^{deg(\phi)}})_{m_i(\lambda_{\phi})}}.\]
The formulas given for conjugacy class size in \cite{Kun} and
\cite{St1} are written in different form; for the reader's benefit
they have been expressed here in the form most useful to us. It
follows that \[ 1+\sum_{n=1}^{\infty} Z_{GL(n,q)} u^n = \prod_{\phi
\neq z} \left[1+ \sum_{n \geq 1} \sum_{\lambda \vdash n}
x_{\phi,\lambda} \frac{u^{n \cdot deg(\phi)}}{q^{deg(\phi) \cdot
\sum_i (\lambda_i')^2} \prod_{i \geq 1}
\left(\frac{1}{q^{deg(\phi)}}\right)_{m_i(\lambda_{\phi})}}\right]. \]
This is called the cycle index generating function.

	Let $Mat(n,q)$ be the set of all $n \times n$ matrices over
the field $F_q$. Define \[ Z_{Mat(n,q)} = \frac{1}{|GL(n,q)|}
\sum_{\alpha \in Mat(n,q)} \prod_{\phi: |\lambda_{\phi}(\alpha)|>0}
x_{\phi,\lambda_{\phi}(\alpha)}. \] Analogous arguments \cite{St1}
show that \[1+\sum_{n=1}^{\infty} Z_{Mat(n,q)} u^n = \prod_{\phi}
\left[1+ \sum_{n \geq 1} \sum_{\lambda \vdash n} x_{\phi,\lambda}
\frac{u^{n \cdot deg(\phi)}}{ q^{deg(\phi) \cdot \sum_i
(\lambda_i')^2} \prod_{i \geq 1}
\left(\frac{1}{q^{deg(\phi)}}\right)_{m_i(\lambda_{\phi})}}\right]. \]
This will be used in Subsection \ref{applic}. Note that the
denominator in $Z_{Mat(n,q)}$ is $|GL(n,q)|$, not $|Mat(n,q)|$, since
the formula follows from a formula for the size of the orbits of
$GL(n,q)$ acting on $Mat(n,q)$ by conjugation. This makes no essential
difference for applications.

\subsection{Applications} \label{applic}

	This subsection describes applications of cycle indices. The
first example is treated in detail and results for the other examples
are sketched.

{\bf Example 1: Cyclic and Separable Matrices}

	Recall that a matrix $\alpha \in Mat(n,q)$ operating on a
vector space $V$ is called cyclic if there is a vector $v_0 \in V$
such that $v_0,v_0 \alpha, v_0 \alpha^2, \cdots$ span $V$. As is
explained in \cite{NP2}, this is equivalent to the condition that the
characteristic and minimal polynomials of $\alpha$ are equal.

	The need to estimate the proportion of cyclic matrices arose
from \cite{NP1} in connection with analyzing the running time of an
algorithm for deciding whether or not the group generated by a given
set of matrices in $GL(n,q)$ contains the special linear group
$SL(n,q)$. Cyclic matrices also arise in recent efforts to improve
upon the MeatAxe algorithm for computing modular characters \cite{NP4}
and in Example 8 below. John Thompson has asked if every matrix is the
product of a cyclic matrix and a permutation matrix, suggesting that
the answer could have applications to finite projective planes.

	Letting $c_M(n,q)$ be the proportion of cyclic elements of
$Mat(n,q)$, the paper \cite{NP2} proves that

\[ \frac{1}{q^2(q+1)} < 1-c_M(n,q) < \frac{1}{(q^2-1)(q-1)}.\] The
cycle index approach is also informative, yielding a formula for the
$n \rightarrow \infty$ limit of $C_M(n,q)$, denoted by
$c_M(\infty,q)$, together with convergence rates. For the argument two
lemmas are useful, as is some notation. Let $N_d(q)$ be the number of
degree $d$ irreducible polynomials over the field $F_q$. In all that
follows $\phi$ will denote a monic irreducible polynomial over
$F_q$. Given a power series $f(u)$, let $[u^n] f(u)$ denote the
coefficient of $u^n$ in $f(u)$.

\begin{lemma}\label{allpoly} \[ \prod_{\phi}
(1-\frac{u^{deg(\phi)}}{q^{deg(\phi)}}) = 1-u \] \end{lemma}

\begin{proof} Expanding $\frac{1}{1-\frac{u^{deg(\phi)}}
{q^{deg(\phi)}}}$ as a geometric series and using unique factorization
in $F_q[x]$, one sees that the coefficient of $u^d$ in the reciprocal
of the left hand side is $\frac{1}{q^d}$ times the number of monic
polynomials of degree $d$, hence 1. Comparing with the reciprocal of
the right hand side completes the proof. \end{proof}

\begin{lemma} \label{bign} If the Taylor series of $f$ around 0
converges at $u=1$, then

\[ lim_{n \rightarrow \infty} [u^n] \frac{f(u)}{1-u} = f(1). \]
\end{lemma}

\begin{proof} Write the Taylor expansion $f(u) = \sum_{n=0}^{\infty}
a_n u^n$. Then observe that $[u^n] \frac{f(u)}{1-u} = \sum_{i=0}^n
a_i$.  \end{proof}

	Theorem \ref{limiting} calculates $c_M(\infty,q)$.

\begin{theorem} \label{limiting} (\cite{F1},\cite{W2}) \[ c_M(\infty,q) =
(1-\frac{1}{q^5}) \prod_{r=3}^{\infty} (1-\frac{1}{q^r})\]
\end{theorem}

\begin{proof} Recall that $\alpha$ is cyclic precisely when its
characteristic polynomial and minimal polynomials are equal. From
Subsection \ref{subsecGL}, these polynomials are equal when all
$\lambda_{\phi}$ have at most one part. In the
cycle index for $Mat(n,q)$ set $x_{\phi,\lambda}=1$ if $\lambda$ has
at most 1 part and $x_{\phi,\lambda}=0$ otherwise. It follows that

\[ c_M(n,q) = \frac{|GL(n,q)|}{q^{n^2}} [u^n] \prod_{\phi}
(1+\sum_{j=1}^{\infty} \frac{u^{j \cdot deg(\phi)}}{q^{(j-1)
deg(\phi)} (q^{deg(\phi)}-1)}).\] By Lemma \ref{allpoly} this equation
can be rewritten as

\begin{eqnarray*}
c_M(n,q) & = &   \frac{|GL(n,q)|}{q^{n^2}} [u^n] \frac{\prod_{\phi} (1-\frac{u^{deg(\phi)}}{q^{deg(\phi)}})
(1+\sum_{j=1}^{\infty} \frac{u^{j \cdot deg(\phi)}}{q^{(j-1)
deg(\phi)} (q^{deg(\phi)}-1)})}{1-u}\\
& = & \frac{|GL(n,q)|}{q^{n^2}} [u^n] \frac{ \prod_{\phi} (1+\frac{u^{deg(\phi)}}{q^{deg(\phi)}(q^{deg(\phi)}-1)})}{1-u}\\
& = & \frac{|GL(n,q)|}{q^{n^2}} [u^n] \frac{ \prod_{d \geq 1} (1+\frac{u^d}{q^{d}(q^d-1)})^{N_d(q)}}{1-u}. 
\end{eqnarray*}

	Recall that a product $\prod_{n=1}^{\infty} (1+a_n)$ converges
absolutely if the series $\sum_{n \geq 1} |a_n|$ converges. Thus using
the crude bound $N_d(q) \leq q^d$ \[ \prod_{d \geq 1}
(1+\frac{u^d}{q^{d}(q^d-1)})^{N_d(q)} \] is analytic in a disc of
radius greater than $1$. Lemma \ref{bign} implies that

\begin{eqnarray*}
c_M(\infty,q) & = & lim_{n \rightarrow \infty} \frac{|GL(n,q)|}{q^{n^2}} [u^n] \frac{ \prod_{d \geq 1} (1+\frac{u^d}{q^{d}(q^d-1)})^{N_d(q)}}{1-u}\\
& = & \prod_{r=1}^{\infty} (1-\frac{1}{q^r}) \prod_{d \geq 1} (1+ \frac{1}{q^{d}(q^{d}-1)})^{N_d(q)}.
\end{eqnarray*} Applying Lemma \ref{allpoly} (with $u=\frac{1}{q}$, $u=\frac{1}{q^2}$ and then $u=\frac{1}{q^5}$) gives that

\begin{eqnarray*}
c_M(\infty,q) & = & \prod_{r=3}^{\infty} (1-\frac{1}{q^r}) \prod_{d \geq 1} ((1+
\frac{1}{q^{d}(q^{d}-1)})(1-\frac{1}{q^{2d}})(1-\frac{1}{q^{3d}}))^{N_d(q)}\\
& = &  \prod_{r=3}^{\infty} (1-\frac{1}{q^r}) \prod_{d \geq 1}
(1-\frac{1}{q^{6d}})^{N_d(q)}\\
& = &  (1-\frac{1}{q^5}) \prod_{r=3}^{\infty} (1-\frac{1}{q^r}).
\end{eqnarray*}
\end{proof}

	The next challenge is to bound the convergence rate of
$c_M(n,q)$ to $c_{\infty}(n,q)$. Wall \cite{W2} found a strikingly
simple way of doing this by relating the cycle index of cyclic
matrices to the cycle index of the set of matrices whose
characteristic polynomial is squarefree (these matrices are termed
separable in \cite{NP2}). To state the result, let $s_M(n,q)$ be the
probability that an $n \times n$ matrix is separable. Next let $C_M(u,q)$ and
$S_M(u,q)$ be the generating functions defined as

\[ C_M(u,q) = 1 + \sum_{n \geq 1} \frac{u^n q^{n^2}}{|GL(n,q)|}
c_M(n,q) \]

\[ S_M(u,q) = 1 + \sum_{n \geq 1} \frac{u^n q^{n^2}}{|GL(n,q)|}
s_M(n,q).\]

\begin{lemma} \label{compare} (\cite{W2}) \[ (1-u) C_M(u,q) =
S_M(u/q,q).\] \end{lemma}

\begin{proof} The proof of Theorem \ref{limiting} shows that

\[ (1-u) C_M(u,q) = \prod_{d \geq 1} (1+\frac{u^{d}}
{q^{d}(q^{d}-1)})^{N_d(q)}.\] A matrix is separable if and only if all
$\lambda_{\phi}$ have size $0$ or $1$. Hence \[ S_M(u,q) =
\prod_{d \geq 1} (1+\frac{u^d}{q^d-1})^{N_d(q)}.\] The result
follows. \end{proof}

\begin{cor} (\cite{W2}) \[ 0 < |c_M(n,q)-c_M(\infty,q)| <
\frac{1}{q^{n+1}(1-1/q)} \] \end{cor}

\begin{proof} Taking coefficients of $u^{n+1}$ on both sides of Lemma
\ref{compare} gives the relation

\[ c_M(n+1,q) - c_M(n,q) = \frac{s_M(n+1,q)-c_M(n,q)}{q^{n+1}}.\]
Since $0 \leq |s_M(n+1,q)-c_M(n,q)| \leq 1$ for all $n$, it follows that

\[|c_M(n,q)-c_M(\infty,q)| \leq \sum_{i=n}^{\infty}
|c_M(i+1,q)-c_M(i,q)| \leq \sum_{i=n}^{\infty} \frac{1}{q^{i+1}},\] as
desired. \end{proof}

{\bf Remarks:}

\begin{enumerate}

\item As mentioned in the introduction, an argument similar to that of
Theorem \ref{limiting} shows that the $n \rightarrow \infty$
probability that an element of $GL(n,q)$ is cyclic is
$(1-\frac{1}{q^5})/(1+\frac{1}{q^3})$. For large $q$ this goes like
$1-1/q^3$. The reason for this is a result of Steinberg \cite{Stei}
stating that the set of non-regular elements in an algebraic group has
co-dimension 3. In type $A$, regular (i.e. centralizer of minimum
dimension) and cyclic elements coincide, but not always. For more
discussion on this point, see \cite{NP2}, \cite{FNP1}.

\item The generating functions $S_M(u,q)$ and $C_M(u,q)$ have
intriguing analytical properties. It is proved in \cite{W2} that

\[ S_M(u,q) = \frac{\prod_{d=1}^{\infty} (1-\frac{u^d(u^d-1)}
{q^d(q^d-1)})^{N_d(q)}}{1-u}.\] Thus $S_M(u,q)$ has a pole at 1 and
$S_M(u,q) - \frac{1}{1-u}$ can be analytically extended to the circle
of radius $q$. Analogous properties hold for $C_M(u,q)$ by means of
Lemma \ref{compare}.

\item The limits $s_M(\infty,q)$ and $s_{GL}(\infty,q)$ are in
\cite{F1},\cite{W2}. Bounding the rate of convergence of $s_M(n,q)$ to
$s_M(\infty,q)$ leads to interesting number theory. Let $p(d)$ be the
number of partitions of $d$ and let $p_2(d)=\sum_{i=0}^d p(d)$. It is
proved in \cite{W2} that \[ |\frac{s_M(n,q) q^{n^2}}{|GL(n,q)|}-1|
\leq \sum_{d=n+1}^{\infty} (p_2(d)+qp(d-2))q^{-d} \leq \frac{1}{3}
(\frac{4q+27}{2q-3}) (\frac{2}{3}q)^{-n}.\]

\item Lehrer \cite{L1} expresses $s_{M}(n,q)$ and $s_{GL}(n,q)$ as
inner products of characters in the symmetric group and proves a
stability result about their expansions in power of $q^{-1}$. See also
\cite{W2} and \cite{LSe}.

\item The results of \cite{F2} and \cite{W2} surveyed above are
extended to the finite classical groups in \cite{FNP1}. The paper
\cite{FlJ} gives (intractable) formulas for the chance of being
separable in groups such as $SL(n,q)$ (i.e. semisimple and simply
connected).

\end{enumerate}

{\bf Example 2: Eigenvalue free matrices}

	The paper \cite{NP3} studies eigenvalues free matrices
(i.e. matrices without fixed lines) over finite fields as a step in
obtaining estimates of cyclic probabilities in orthogonal groups
\cite{NP5}. It is interesting that the study of eigenvalue free
matrices was one of the motivations for the original papers
\cite{Kun},\cite{St1}, the latter of which proves that the $n,q
\rightarrow \infty$ limit of the chance that an element of $GL(n,q)$
has no eigenvalues is $\frac{1}{e}$.

	The $n \rightarrow \infty$ probability that a random element
of $S_n$ has no fixed points is also $\frac{1}{e}$. This is not
coincidence; in general the $q \rightarrow \infty$ limit of the chance
that the characteristic polynomial of a random element of $Mat(n,q)$
factors into $n_i$ degree $i$ irreducible factors is the same as the
probability that an element of $S_n$ factors into $n_i$ cycles of
degree $i$. This is proved at the end of \cite{St1} and is extended to
finite Lie groups in \cite{F1} using the combinatorics of maximal
tori. There is another interesting line of argument which should be
mentioned. It is easy to see from the cycle index that the
factorization type of the characteristic polynomial of a random
element of $Mat(n,q)$ and the factorization type of a random degree
$n$ polynomial over $F_q$ have the same distribution as $q \rightarrow
\infty$. Now the factorization type of a random degree $n$ polynomial
over $F_q$ has same distribution as the cycle type of a random
permutation distributed as a $q$-shuffle on $n$ cards \cite{DMP}, and
as $q \rightarrow \infty$ a $q$-shuffle converges to a random
permutation. The connection of Lie theory with card shuffling may seem
adhoc, but is really the tip of a deep iceberg \cite{F10}.

{\bf Example 3: Characteristic polynomials}

	The previous example is a special case of the problem of
studying the degrees of the factors of the characteristic polynomial
of a random matrix. Many results in this direction (all proved used
cycle indices) can be found in Stong's paper \cite{St1}. Hansen and
Schmutz \cite{HSc} use cycle index manipulations to prove that if one
ignores factors of small degree, then the factorization type of the
characteristic polynomial of a random element of $GL(n,q)$ is close to
the factorization type of a random degree $n$ polynomial over
$F_q$. More precisely, let $A_{n,l}$ be the set of sequences
$(\alpha_{l+1},\cdots,\alpha_n)$ where $\alpha_i$ is the number of
degree $i$ factors of a random polynomial chosen from some
measure. Let $Q_n^{(1)}$ be the measure on polynomials arising from
characteristic polynomials of random elements of $GL(n,q)$ and let
$Q_n^{(2)}$ be the measure arising from choosing a degree $n$
polynomial over $F_q$ uniformly at random. They prove

\begin{theorem} (\cite{HSc}) There exists constants $c_1,c_2$ such that
for all $l$ with $c_1log(n) \leq l \leq n$ and $B \subset N^{n-l}$,

\[ |Q_n^{(1)}(A_n(B))-Q_n^{(2)}(A_n(B))| < c_2/l. \] \end{theorem} The
final section of their paper uses this principle to prove results
about characteristic polynomials of random matrices using known
results about random polynomials. A useful reference on the
distribution of degrees of random polynomials over finite fields is
\cite{ABT}.

{\bf Example 4: Generating transvections}

	Recall that the motivation behind Example 1 was a group
recognition problem, i.e. trying to determine whether or not the group
generated by a given set $X$ of matrices in $GL(n,q)$ contains the
special general linear group $SL(n,q)$. However the problem still
remains of making the recognition algorithm constructive. For instance
if the group generated by $X$ is $GL(n,q)$ it would be desirable to
write any element of $GL(n,q)$ as a word in $X$.

	The paper \cite{CL} proposes such a constructive recognition
algorithm. An essential step involves constructing a transvection,
that is a non-identity element of $SL(n,q)$ which has an $n-1$
dimensional fixed space. This in turn is done in two steps. First,
find an element $\alpha$ of $GL(n,q)$ conjugate to
diag$(C((z-\tau)^2),R)$ where $C$ is the companion matrix as in
Subsection \ref{subsecGL} and $R$ is semisimple without $\tau$ as an
eigenvalue. Second, one checks that raising $\alpha$ to the least
common multiple of the orders of $\tau$ and $R$ gives a transvection.

	Thus it necessary to bound the number of feasible $\alpha$ in
the first step. Such $\alpha$ have conjugacy class data
$\lambda_{z-\tau}=(2)$ and all other $\lambda_{\phi}$ have largest
part at most 1. The cycle index approach gives bounds improving on
those in \cite{CL}; see \cite{FNP1} for the details.

{\bf Example 5: Semisimple matrices} 

	A fundamental problem in computational group theory is to
construct an element of order $p$. Given a group element $g$ with
order a multiple of $p$, this can be done by raising $g$ to an
appropriate power. It is proved in \cite{IKS} that if $G$ is a
permutation group of degree $n$ with order divisible by $p$, then the
probability that a random element of $G$ has order divisible by $p$ is
at least $\frac{1}{n}$.

	Their proof reduces the assertion to simple groups and then
uses the classification of simple groups. Let us consider the group
$GL(n,q)$, which is close enough to simple to be useful for the
applications at hand. When $p$ is the characteristic of the field of
definition of $GL(n,q)$, an element has order prime to $p$ precisely
when it is semisimple. Thus the problem is to study the probability
that an element of $GL(n,q)$ is semisimple. The paper \cite{GL} shows
that if $G$ is a simple Chevalley group, then the probability of not
being semisimple is at most $3/(q-1)+2/(q-1)^2$ and thus at most $c/q$
for some constant $c$ as conjectured by Kantor.

	As mentioned earlier, a matrix $\alpha$ is semisimple if and
only if all $\lambda_{\phi}(\alpha)$ have largest part size at most
1. Stong \cite{St1} used cycle indices to obtain crude asymptotic
bounds for the probability that an element of $GL(n,q)$ is
semisimple. The thesis \cite{F1} used the Rogers-Ramanujan identities
to prove that the $n \rightarrow \infty$ probability that an element
of $GL(n,q)$ is semisimple is

\[ \prod_{r=1 \atop r=0,\pm 2 (mod \ 5)}^{\infty}
\frac{(1-\frac{1}{q^{r-1}})}{(1-\frac{1}{q^r})}. \] The paper
\cite{FNP1} gives effective bounds for finite $n$.

{\bf Example 6: Order of a matrix}

	A natural problem is to study the order of a random
matrix. This has been done in \cite{St2} and \cite{Schm}; see also the
remarks in Subsection \ref{samplingsize} and the very preliminary
calculations for other classical groups in \cite{F1}. Shalev
\cite{Sh1} uses facts about the distribution of the order of a random
matrix together with Aschbacher's study of maximal subgroups of
classical groups \cite{As} as key tools in studying the probability
that a random element of $GL(n,q)$ belongs to an irreducible subgroup
of $GL(n,q)$ that does not contain $SL(n,q)$. As explained in
\cite{Sh1} this has a number of appications; for instance it leads to
a proof that if $x$ is any non-trivial element of $PSL(n,q)$ then the
probability that $x$ and a randomly chosen element $y$ generate
$PSL(n,q)$ tends to $1$ as $q \rightarrow \infty$. Shalev \cite{Sh1}
asks for extensions of these results to other finite classical groups.

	It is also useful to count elements of given orders (e.g. $2$
or $3$) in classical groups and their maximal subgroups. The recent
paper \cite{CTY} uses cycle indices to perform such enumerations. One
motivation for such enumerations is the study of finite simple
quotients of $PSL(2,Z)$; a group $G$ is a quotient of $PSL(2,Z)$ if
and only if $G=<x,y>$ with $x^2=y^3=1$. For further discussion, see
\cite{Sh2}.

{\bf Example 7: Random number generators}

	We follow \cite{Ma},\cite{MaT} in indicating the relevance of
random matrix theory to the study of random number generators. Suppose
one wants to test a mechanism for generating a random integer between
$0$ and $2^{33}-1$. In base $2$ these are length $33$ binary
vectors. Generating say $n$ of these and listing them gives an $n
\times 33$ matrix. If the random generator were perfect, the arising
matrix would be random. One could choose a statistic such as the rank
of a matrix and compare the generation method with theory. They report
that shift-register generators will fail such tests but that
congruential generators usually pass. It would be interesting to see
how various random number generators perform when tested using other
conjugacy class functions of random matrices.

	Diaconis and Graham \cite{DG} analyze random walks of the form
$X_n = A X_{n-1} + \epsilon_n$ where $X_i$ is a length $d$ $0-1$
vector, $A$ is an element of $GL(n,2)$, and $\epsilon_n$ is a random
vector of disturbance terms. For more general $A$ (in $GL(n,q)$) this
includes the problem of running a psuedo-random number generator with
recurrence $Y_n=a_1 Y_{n-1} + \cdots + a_d Y_{n-d} + \epsilon_{n}$
with $Y_i \in F_q$. They show that the rational canonical form of $A$
is related in a subtle way to the convergence rate of the walk. It
would be interesting to understand what happens when $A$ is a random
matrix.

{\bf Example 8: Product replacement algorithm}

	In recent years finite group theory has become much more
computational. Given a generating set $S$ of a finite group $G$, it
is natural to seek random elements of $G$. One approach, implemented
in the computer systems GAP and MAGMA, is the product replacement
algorithm \cite{CeLgMuNiOb}. Fixing $G$ and some $k$, one performs a
random walk on $k$-tuples $(g_1,\cdots,g_k)$ of elements of $G$ which
generate the group. The walk proceeds by picking an ordered pair
$(i,j)$ with $1 \leq i \neq j \leq n$ uniformly at random and applying one
of the following four operations with equal probability:

\[ R_{i,j}^{\pm}: (g_1,\cdots,g_i,\cdots,g_k) \mapsto (g_1,\cdots,g_i
\cdot g_j^{\pm},\cdots,g_k) \]

\[ L_{i,j}^{\pm}: (g_1,\cdots,g_i,\cdots,g_k) \mapsto (g_1,\cdots,
g_j^{\pm} \cdot g_i,\cdots,g_k). \] These moves map generating
$k$-tuples to generating $k$-tuples. One starts from any generating
$k$-tuple, applies the algorithm for $r$ steps, and then outputs a
random entry of the resulting $k$-tuple (i.e. a group element).

	The product replacement algorithm has superb practical
performance (often converging more rapidly than random walk on the
Cayley graph), in spite of the theoretical defects that a random entry
of a random generating $k$-tuple does not have the same distribution
as a random element of $G$, and that the convergence rate of the chain
on $k$-tuples to its stationary distribution is unknown. The paper
\cite{CeLgMuNiOb}, aware of these issues, tests the algorithm against
theory, using conjugacy class statistics such as the order of an
element, the number of factors of the characteristic polynomial of a
random matrix, the degree of the largest irreducible factor of the
characteristic polynomial of a random matrix, and the proportion of
cyclic matrices in the finite classical groups. In short,
understanding properties of random matrices is crucial to their
analysis.

	A recent effort to understand the performance of the product
replacement algorithm uses Kazhdan's property (T) from the
representation theory of Lie groups \cite{LP}; the paper \cite{P} is a
useful survey. Much remains to be done.

{\bf Example 9: Running times of algorithms} 

	One of the main approaches to computing determinants and
permanents of integer matrices involves doing the computations for
reductions mod prime powers. Section 4.6.4 of \cite{Kn} gives a
detailed discussion with references to literature on upper bounds of
running times. If one believes that typical matrices one encounters in
the real world are like random matrices, this motivates studying
random matrices over finite fields. In fact von Neumann's interest in
eigenvalues of random matrices with independent normal entries arose
from the same heuristic applied to questions in numerical analysis
(the introduction of \cite{E} gives further discussion of this point).

	Examples of algorithms in which properties of random matrices
were really needed to bound running times include recognizing when a
group generated by a set of matrices contains $SL(n,q)$ \cite{NP1} and
the MeatAxe algorithm for computing modular characters \cite{NP4}.

{\bf Example 10: Isometry classes of linear codes}

	Fripertinger \cite{Frip1}, \cite{Frip2} considers cycle
indices (in the permutation sense) of matrix groups acting on
lines. His interest was in understanding properties of random isometry
classes of linear codes--a harder problem than understanding random
linear codes. The cycle indices he obtains seem quite intractable for
theorem proving, but are useful in conjunction with computers. He also
gives references to the switching function literature.

	Curiously, understanding the permutation action of random
matrices of lines comes up in another context. Wieand \cite{W} has
shown that the eigenvalues of random permutation matrices possess a
structure similar to the eigenvalues of matrices from compact Lie
groups. Persi Diaconis has suggested that the eigenvalues of high
dimensional representations of finite groups of Lie type (such as the
permutation action on lines) may possess similar structure; see
\cite{F5} for more in this direction.
	
\subsection{Generalization to the Classical Groups} \label{Classical}

	This subsection will focus on the finite unitary groups, with
remarks about symplectic and orthogonal groups at the end. These cycle
indices were derived in \cite{F1},\cite{F2} and were applied to the
problem of estimating proportions of cyclic, separable, and semisimple
matrices (these terms were defined in Subsection \ref{applic}) in
\cite{FNP1}. 
	
	The unitary group $U(n,q)$ can be defined as the subgroup of
$GL(n,q^2)$ preserving a non-degenerate skew-linear form. Recall that
a skew-linear form on an $n$ dimensional vector space $V$ over
$F_{q^2}$ is a bilinear map $<,>:V \times V \rightarrow F_{q^2}$ such
that $<\vec{x},\vec{y}> = <\vec{y},\vec{x}>^q$ (raising to the $q$th
power is an involution in a field of order $q^2$). One such form is
given by $<\vec{x},\vec{y}> = \sum_{i=1}^n x_i y_i^q$. Any two
non-degenerate skew-linear forms are equivalent, so that $U(n,q)$ is
unique up to isomorphism.

	Wall $\cite{W1}$ parametrized the conjugacy classes of the
finite unitary groups and computed their sizes. To describe his
result, an involution on polynomials with non-zero constant term is
needed. Given a polynomial $\phi$ with coefficients in $F_{q^2}$ and
non vanishing constant term, define a polynomial $\tilde{\phi}$ by:

\[ \tilde{\phi} = \frac{z^{deg(\phi)}
\phi^q(\frac{1}{z})}{[\phi(0)]^q} \] where $\phi^q$ raises each
coefficient of $\phi$ to the $q$th power. Writing this out, a
polynomial $\phi(z)=z^{deg(\phi)} + \alpha_{deg(\phi)-1}
z^{deg(\phi)-1} + \cdots + \alpha_1 z + \alpha_0$ with $\alpha_0 \neq
0$ is sent to $\tilde{\phi}(z)= z^{deg(\phi)} +
(\frac{\alpha_1}{\alpha_0})^q z^{deg(\phi)-1}+ \cdots +
(\frac{\alpha_{deg(\phi)-1}} {\alpha_0})^qz +
(\frac{1}{\alpha_0})^q$. An element $\alpha \in U(n,q)$ associates to
each monic, non-constant, irreducible polynomial $\phi$ over $F_{q^2}$
a partition $\lambda_{\phi}$ of some non-negative integer
$|\lambda_{\phi}|$ by means of rational canonical form. The
restrictions necessary for the data $\lambda_{\phi}$ to represent a
conjugacy class are that $|\lambda_z|=0$,
$\lambda_{\phi}=\lambda_{\tilde{\phi}}$, and that $\sum_{\phi}
|\lambda_{\phi}|deg(\phi)=n.$

	Using formulas for conjugacy class sizes from \cite{W1}
together with some combinatorial manipulations, the following unitary
group cycle index generating function was derived in \cite{F1}. The
products in the theorem are as always over monic irreducible
polynomials.

\begin{theorem}
\begin{eqnarray*}
& & 1+\sum_{n=1}^{\infty} \frac{u^n}{|U(n,q)|} \sum_{\alpha \in
U(n,q)} \prod_{\phi: |\lambda_{\phi}(\alpha)|>0} x_{\phi,\lambda_{\phi}(\alpha)} \\
& = & \prod_{\phi \neq z,
\phi=\tilde{\phi}} \left[1+ \sum_{n \geq 1} \sum_{\lambda \vdash n} x_{\phi,\lambda}
\frac{(-u)^{n \cdot deg(\phi)}}{ (-q)^{deg(\phi) \cdot \sum_i
(\lambda_i')^2} \prod_{i \geq 1} \left(\frac{1}{(-q)^{deg(\phi)}}\right)_{m_i(\lambda)}}\right]\\
& & \cdot \prod_{\{\phi,\tilde{\phi}\}, \phi \neq \tilde{\phi}} \left[1+\sum_{n \geq 1}\sum_{\lambda \vdash n}
x_{\phi,\lambda} x_{\tilde{\phi},\lambda}
\frac{u^{2n \cdot deg(\phi)}}{ q^{2deg(\phi) \cdot \sum_i
(\lambda_i')^2} \prod_{i \geq 1} \left(\frac{1}{q^{2deg(\phi)}}\right)_{m_i(\lambda)}  }\right]
\end{eqnarray*}
\end{theorem}

	One interesting theoretical result concerning the cycle index
of $U(n,q)$ is the following functional equation. Letting
$C_{GL}(u,q)$ and $C_{U}(u,q)$ be the cycle index generating functions
for cyclic matrices in the general linear and unitary groups
respectively, the functional equation states that \[ C_{GL}(u,q)
C_U(-u,-q) = C_{GL}(u^2,q^2) .\] The paper \cite{FNP1} proves that
this relation holds whenever the condition on the partitions
$\lambda_{\phi}$ is independent of the polynomial $\phi$. In the
current example, a matrix is cyclic if and only if all
$\lambda_{\phi}$ have at most one row. This condition is independent
of $\phi$.

	Cycle indices for the symplectic and orthogonal groups are a
bit trickier to establish from Wall's formulas. To the treatment in
\cite{F1},\cite{F2} we add a remark which should be very helpful to
anyone trying to use those cycle indices. The paper \cite{F2} only
wrote out an explicit formula for the cycle index for the sum of
$+$,$-$ type orthogonal groups. To solve for an individual orthogonal
group, it is necessary to average that formula with a formula for the
difference of $+$,$-$ type orthogonal groups (this procedure is
carried out in a special case in \cite{FNP1}). In general, the formula
for the difference of orthogonal groups is obtained from the formula
for the sum of orthogonal groups as follows. First, for the
polynomials $z \pm 1$, replace terms corresponding to partitions with
an odd number of odd parts by their negatives. Second, for polynomials
invariant under $\tilde{}$, replace terms corresponding to partitions
of odd size by their negatives.

\subsection{Limitations and Other Methods}

	Cycle index techniques, while very useful, also have their
limitations and are not always the best way to proceed, as the
following examples demonstrate.

{\bf Example 1: Primitive prime divisor elements}

	For integers $b,e>1$ a primitive prime divisor of $b^e-1$ is a
prime dividing $b^e-1$ but not dividing $b^i-1$ for any $i$ with $1
\leq i <e$. An element of $GL(n,q)$ is called a primitive prime
divisor (ppd) element if its order is divisible by a primitive prime
divisor of $q^e-1$ with $n/2 < e \leq n$. The analysis in \cite{NiP}
derives elegant bounds on the proportions of ppd elements in the
finite classical groups and applies them to the group recognition
problem for classical groups over finite fields (determining when a
group generated by a set of matrices contains $SL(n,q)$). We do not
see how to get comparable bounds using generating function techniques.

{\bf Example 2: Proportions of semisimple elements in exceptional
groups}

	Although Example 5 of Section \ref{applic} was estimating
proportions of semisimple matrices, this was only for the finite
classical groups, where the index $n$ can take an infinite number of
values. Cycle indices don't seem useful unless there is a tower of
groups of varying rank available.

	Fortunately the computer package CHEVIE permits calculations
precisely in finite rank cases such as the exceptional groups. Indeed
this is how \cite{GL} obtained estimates of the proportions of
semisimple elements in the exceptional groups.

{\bf Example 3: Non-uniform distributions on matrices}

	The cycle indices give useful information about conjugacy
class functions when the matrix is chosen uniformly at random. However
there are other distributions on matrices which one could study and
for which cycle index methods (at present) can not be applied.

	One example is random $n \times n$ matrices where the matrix
entries are chosen independently according to a given probability
distribution on $F_q$. Charlap, Rees, and Robbins \cite{CRR} show that
if the probability distribution is not concentrated on any proper
affine subspace of $F_q$, then as $n \rightarrow \infty$ the
probability that the matrix is invertible is the same as for a uniform
matrix. They use Moebius inversion on the lattice of subspaces of an
$n$ dimensional vector space and the Poisson summation formula. Is the
same true for other natural conjugacy class functions? We expect that
the answer is yes, which can be regarded as a type of ``universality''
result for the asymptotic description of random elements of $GL(n,q)$
to be given in Subsection \ref{measures}. Analogous universality
results are known for matrices with complex entries \cite{So}. For
further information on the rank of random $0-1$ matrices, see
\cite{BKW} for sparse matrices, \cite{Bo} for a survey of results on
the rank over the real numbers, and also the discussion of work of
Rudvalis and Shinoda in Subsection \ref{Symmetric}.

	It is conceivable that cycle index techniques will be able to
handle certain natural non-uniform distributions on $GL(n,q)$. This
happens for the symmetric groups, where natural non-uniform measures
such as performing a $q$-riffle shuffle on a deck of cards has a
useful cycle index \cite{DMP},\cite{F10}.

\section{Running Example: General Linear Groups} \label{GLSec}

	The purpose of this section is to give different ways of
understanding the conjugacy class of a random element of
$GL(n,q)$. The analogous theory for other finite classical groups is
mentioned in passing but is not treated in detail as many of the main
ideas can be communicated using $GL(n,q)$. Subsection \ref{measures}
will show how this leads naturally to the study of certain probability
measures $M_{GL,u,q}$ on the set of all partitions of all natural
numbers. Connections with symmetric function theory lead to several
ways of growing random partitions distributed according to
$M_{GL,u,q}$. One consequence is a motivated proof of the
Rogers-Ramanujan identities.

\subsection{Measures on Partitions} \label{measures}

	The goal is to obtain a probabilistic description of the
conjugacy class of a random element of $GL(n,q)$. The ideas are based
on \cite{F1}. For this the following definition will be fundamental.

{\bf Definition:} The measure $M_{GL,u,q}$ on the set of all
partitions of all natural numbers is defined by

\[ M_{GL,u,q}(\lambda) = \prod_{r=1}^{\infty} (1-\frac{u}{q^r})
\frac{u^{|\lambda|}}{q^{\sum_i (\lambda_i')^2} \prod_i
(\frac{1}{q})_{m_i(\lambda)}}.\]

	The motivation for this definition will be clear from Theorem
\ref{relate}. The measure $M_{GL,u,q}$, while seemingly complicated,
does have some nice combinatorial properties. For instance for
partitions of a fixed size, this measure respects the dominance order
on partitions (in this partial order $\lambda \geq \mu$ if and only if
$\lambda_1 + \cdots \lambda_i \geq \mu_1 + \cdots + \mu_i$ for all $i
\geq 1$). In work with Bob Guralnick we actually needed this
property.

	Lemma \ref{probabmeasure} proves that for $q>1$ and $0<u<1$,
the measure $M_{GL,u,q}$ is in fact a probability measure. There are
at least three other proofs of this fact: an argument using $q$
series, specializing an identity about Hall-Littlewood polynomials, or
a slick argument using Markov chains and an identity of Cauchy. This
third argument will be given in Subsection \ref{MarkovGL}.

\begin{lemma} \label{probabmeasure} If $q>1$ and $0<u<1$, then
$M_{GL,u,q}$ defines a probability measure. \end{lemma}

\begin{proof} $M_{GL,u,q}$ is clearly non-negative when $q>1$ and
$0<u<1$. Stong \cite{St1} established an equation which is equivalent
to the sought identity

\[ \sum_{\lambda}\frac{u^{|\lambda|}} {q^{\sum_i (\lambda_i')^2}
\prod_i (\frac{1}{q})_{m_i(\lambda)}} = \prod_{r=1}^{\infty}
(\frac{1}{1-\frac{u}{q^r}}).\] As some effort is required to see this
equivalence, we derive the identity directly using Stong's line of
reasoning.

	First observe that unipotent elements of $GL(n,q)$
corresponding to nilpotent $n \times n$ matrices (subtract the
identity matrix), and that the number of nilpotent $n \times n$
matrices is $q^{n(n-1)}$ by the Fine-Herstein theorem \cite{FeiH}. The
number of unipotent elements in $GL(n,q)$ can be evaluated in another
way using the cycle index of the general linear groups. Namely set
$x_{\phi,\lambda}=1$ if $\phi=z-1$ and set  $x_{\phi,\lambda}=0$
otherwise. One concludes that \[ \sum_{\lambda \vdash n}\frac{1}
{q^{\sum_i (\lambda_i')^2} \prod_i (\frac{1}{q})_{m_i(\lambda)}} =
\frac{q^{n(n-1)}}{|GL(n,q)|}.\] Now multiply both sides by $u^n$, sum
in $n$, and apply Euler's identity \[ \sum_{n=0}^{\infty}
\frac{u^nq^{{n \choose 2}}}{(q^n-1) \cdots (q-1)} =
\prod_{r=1}^{\infty} (\frac{1}{1-\frac{u}{q^r}}).\] \end{proof}

	The measure $M_{GL,u,q}$ is a fundamental object for
understanding the probability theory of conjugacy classes of
$GL(n,q)$.  This emerges from Theorem \ref{relate}.

\begin{theorem} \label{relate} 

\begin{enumerate}

\item Fix $u$ with $0<u<1$. Then choose a random natural number $N$
with probability of getting $n$ equal to $(1-u)u^n$. Choose $\alpha$
uniformly in $GL(N,q)$. Then as $\phi$ varies, the random partitions
$\lambda_{\phi}(\alpha)$ are independent random variables, with
$\lambda_{\phi}$ distributed according to the measure
$M_{GL,u^{deg(\phi)},q^{deg(\phi)}}$.

\item Choose $\alpha$ uniformly in $GL(n,q)$. Then as $n \rightarrow
\infty$, the random partitions $\lambda_{\phi}(\alpha)$ converge in
finite dimensional distribution to independent random variables, with
$\lambda_{\phi}$ distributed according to the measure
$M_{GL,1,q^{deg(\phi)}}$.

\end{enumerate}
\end{theorem}

\begin{proof} Recall the cycle index factorization

 \[ 1+\sum_{n=1}^{\infty} Z_{GL(n,q)} u^n = \prod_{\phi \neq z}
\left[1+\sum_{n \geq 1} \sum_{\lambda \vdash n} x_{\phi,\lambda}
\frac{u^{n \cdot deg(\phi)}}{\prod_{\phi} q^{deg(\phi) \cdot \sum_i
(\lambda_i')^2} \prod_{i \geq 1} (\frac{1}{q^{deg(\phi)}})_{m_i}}\right]. \] Setting all
$x_{\phi,\lambda}$ equal to $1$ and using Lemma
\ref{probabmeasure} shows that

\[ \frac{1}{1-u} = \prod_{\phi \neq z} \prod_{r=1}^{\infty}
(\frac{1}{1-\frac{u^{deg(\phi)}}{q^{r \cdot deg(\phi)}}}).\] Taking reciprocals
and multiplying by the cycle index factorization shows that

\[ (1-u)+\sum_{n=1}^{\infty} Z_{GL(n,q)} (1-u) u^n = \prod_{\phi \neq
z} \left( M_{GL,u^{deg(\phi)},q^{deg(\phi)}}(\emptyset)+\sum_{\lambda: |\lambda|>0}
M_{GL,u^{deg(\phi)},q^{deg(\phi)}}(\lambda) x_{\phi,\lambda}\right).\] This
proves the first assertion of the theorem. For the second assertion,
use Lemma \ref{bign} from Subsection \ref{applic}.  \end{proof}

{\bf Remarks:}

\begin{enumerate}

\item Theorem \ref{relate} has an analog for the symmetric groups
\cite{SL}. The statement is as follows. {\it Fix $u$ with
$0<u<1$. Then choose a random natural number $N$ with probability of
getting $n$ equal to $(1-u)u^n$. Choose $\pi$ uniformly in
$S_N$. Letting $n_i$ be the number of $i$-cycles of $\pi$, the random
variables $n_i$ are independent, with $n_i$ distributed as a Poisson
with mean $\frac{u^i}{i}$. Furthermore if one chooses $\pi$ uniformly
in $S_n$ and lets $n \rightarrow \infty$, then the random variables
$n_i$ are independent random variables, with $n_i$ distributed as a
Poisson$(\frac{1}{i})$.}

\item The idea of performing an auxiliary randomization of $n$ is a
mainstay of statistical mechanics, known as the grand canonical
ensemble. For a clear discussion see Sections 1.7, 1.9, and 4.3 of
\cite{Fe}.

\end{enumerate}

\subsection{Symmetric Function Theory and Sampling Algorithms}
\label{Symmetric}

	The aim of this subsection is two-fold. First, the measures
$M_{GL,u,q}$ are connected with the Hall-Littlewood symmetric
functions. Then we indicate how this connection can be exploited to
give probabilistic methods for growing random partitions distributed
as $M_{GL,u,q}$. The purpose is not to drown the reader in formulas,
but rather to show that the connection between symmetric functions and
probability is deep, beautiful, and useful in both directions. The
results on this section are based on \cite{F1} and \cite{F3}, except
for the remark on how to make the algorithms terminate in finite time,
which is joint with Mark Huber.
 
	To begin, we recall the Hall-Littlewood symmetric functions,
which arise in many parts of mathematics: enumeration of $p$ groups,
representation theory of $GL(n,q)$, and counting automorphisms of
modules. The basic references for Hall-Littlewood polynomials
$P_{\lambda}$ is Chapter 3 of \cite{Mac}, which offers the following
definition \[ P_{\lambda}(x_1,\cdots,x_n;t) = \left[\frac{1}{\prod_{i
\geq 0} \prod_{r=1}^{m_i(\lambda)} \frac{1-t^r}{1-t}}\right] \sum_{w
\in S_n} w \left(x_1^{\lambda_1} \cdots x_n^{\lambda_n} \prod_{i<j}
\frac{x_i-tx_j} {x_i-x_j}\right). \] Here $w$ is a permutation acting
on the $x$-variables by sending $x_i$ to $x_{w(i)}$. Recall that
$m_i(\lambda)$ is the number of parts of $\lambda$ of size $i$. At
first glance it is not obvious that these are polynomials, but the
denominators cancel out after the symmetrization. The Hall-Littlewood
polynomials interpolate between the Schur functions ($t=0$) and the
monomial symmetric functions ($t=1$).

	Theorem \ref{relateHL} relates the measures $M_{GL,u,q}$ to
the Hall-Littlewood polynomials. Recall that $n(\lambda) = \sum_{i}
(i-1) \lambda_i = \sum_{i} {\lambda_i' \choose 2}$.

\begin{theorem} \label{relateHL}

\[ M_{GL,u,q}(\lambda) = \prod_{i=1}^{\infty} (1-\frac{u}{q^{i}})
\frac{P_{\lambda}(\frac{u}{q},\frac{u}{q^{2}} ,
\cdots;\frac{1}{q})}{q^{n(\lambda)}} \] \end{theorem}

\begin{proof} From the above formula for Hall-Littlewood polynomials,
it is clear that the only surviving term in the specialization
$P_{\lambda}(\frac{u}{q},\frac{u}{q^{2}} , \cdots;\frac{1}{q})$ is the
term when $w$ is the identity. The rest is a simple combinatorial
verification. (Alternatively, one could use ``principal
specialization'' formulas for Macdonald polynomials on page 337 of
\cite{Mac}). \end{proof}

{\bf Remark:} The paper \cite{F3} gives symmetric function theoretic
generalizations of the measure $M_{GL,u,q}$ on partitions. In the case
of Schur functions $s_{\lambda}$, this measure depends on two infinite
sets of variables $x_i,y_i$ and assigns a partition $\lambda$ mass
equal to $ s_{\lambda}(x_i) s_{\lambda}(y_i) \prod_{i,j}
(1-x_iy_j)$. It is remarkable that precisely this measure arose in
work of the random matrix community relating the distribution of the
lengths of increasing subsequences of random permutations to the
distribution of eigenvalues of random GUE matrices (these matrices
have complex entries). To elaborate, the Robinson-Schensted-Knuth
correspondence associates a random partition of size $n$ to a random
permutation of size $n$ and the shape of the partition encodes
information about the longest increasing subsequence of the
permutation. Choosing the size of the symmetric group randomly
(according to a Poisson distribution) gives a probability measure on
the set of all partitions of all natural numbers which is a special
case of the above Schur function measure. Then the coordinate change
$h_j=\lambda_1'+\lambda_j-j$ maps the set of row lengths
$\{\lambda_j\}$ of the partition to a set of distinct integers
$\{h_j\}$. These $h_j$ can be viewed as positions of electrostatic
charges repelling each other, and from this viewpoint the measure on
subsets of the integers bears a striking resemblance to the eigenvalue
density of a random GUE matrix. This fantastic heuristic can be made
precise and led to a solution of the long-standing conjecture relating
lengths of increasing subsequences of permutations to eigenvalues of
random matrices. For these developments see \cite{BOO},\cite{Jo} and
the many references therein.

	Now we return to the measure $M_{GL,u,q}$ and describe an
algorithm for growing random partitions according to this measure.

\begin{center}
The Young Tableau Algorithm
\end{center}

\begin{description}

\item [Step 0] Start with $N=1$ and $\lambda$ the empty
partition. Also start with a collection of coins indexed by the
natural numbers, such that coin $i$ has probability $\frac{u}{q^i}$ of
heads and probability $1-\frac{u}{q^i}$ of tails.

\item [Step 1] Flip coin $N$.

\item [Step 2a] If coin $N$ comes up tails, leave $\lambda$ unchanged,
set $N=N+1$ and go to Step 1.

\item [Step 2b] If coin $N$ comes up heads, choose an integer $S>0$
according to the following rule. Set $S=1$ with probability $\frac
{q^{N-\lambda_1'}-1} {q^N-1}$. Set $S=s>1$ with probability
$\frac{q^{N-\lambda_s'}-q^{N-\lambda_{s-1}'}}{q^N-1}$. Then increase
the size of column $s$ of $\lambda$ by 1 and go to Step 1.

\end{description}

        As an example of the Young Tableau Algorithm, suppose we are
at Step 1 with $\lambda$ equal to the following partition:
        
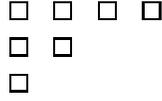
\[ \begin{array}{c c c c}
                \framebox{}& \framebox{}& \framebox{}& \framebox{}  \\
                \framebox{}& \framebox{}&&      \\
                \framebox{} &&&
          \end{array} \] Suppose also that $N=4$ and that coin 4 had
already come up heads once, at which time we added to column 1, giving
$\lambda$. We flip coin 4 again and get heads, going to Step 2b. We
add a box to column $1$ with probability $\frac{q-1}{q^4-1}$, to
column $2$ with probability $\frac{q^2-q}{q^4-1}$, to column $3$ with
probability $\frac{q^3-q^2} {q^4-1}$, to column $4$ with probability
$0$, and to column $5$ with probability $\frac{q^4-q^3}{q^4-1}$. We
then return to Step 1.

\begin{theorem} \label{TableauAlg} For $0<u<1$ and $q>1$, the Young
Tableau Algorithm generates partitions which are distributed according
to the measure $M_{GL,u,q}$. \end{theorem}

	To give insight into the proof of Theorem \ref{TableauAlg}, we
remark that it was deduced by proving a stronger result (Theorem
\ref{stronger}) inductively and then taking the $N \rightarrow \infty$
limit. As is clear from the statement of Theorem \ref{stronger}, the
connection with Hall-Littlewood polynomials (in particular the ability
to truncate them) was crucial. It is unlikely that the Young Tableau
Algorithm would have been discovered without this connection.

\begin{theorem} \label{stronger} Let $P^N(\lambda)$ be the probability
that the algorithm outputs $\lambda$ when coin $N$ comes up tails. Then

\[ P^N(\lambda) =             \left\{ \begin{array}{ll}
                                                                                                                                                                                                                                                                                               \frac{u^{|\lambda|} (\frac{u}{q})_N
(\frac{1}{q})_N}{(\frac{1}{q})_{N-\lambda_1'}} \frac{P_{\lambda}(\frac{1}{q},\cdots,\frac{1}{q^N},0,\cdots;0,\frac{1}{q})}{q^{n(\lambda)}} & \mbox{if $\lambda_1' \leq N$}\\
                                                                                                                                                                                                                                                                                                0 & \mbox{if $\lambda_1' > N$}.
                                                                                                                                                                                                                                                                                                \end{array}
                        \right.                  \] \end{theorem}

	Next we explain why the Young Tableau Algorithm is called
that. A standard Young tableau $T$ of size $n$ is a partition of $n$
with each box filled by one of $\{1,\cdots,n\}$ such that each of
$\{1,\cdots,n\}$ appears exactly once and the numbers increase in each
row and column of $T$. For instance,

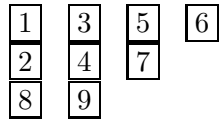
\[ \begin{array}{c c c c c}
                \framebox{1} & \framebox{3} & \framebox{5} & \framebox{6} &   \\
                \framebox{2} & \framebox{4} & \framebox{7} &  &    \\
                \framebox{8} & \framebox{9} &  &  &    
          \end{array} \] is a standard Young tableau. Standard Young
tableaux are important in combinatorics and representation theory. The
Young Tableau Algorithm is so named because numbering the boxes in the
order in which they are created gives a standard Young tableau. Thus
although our initial interest was in the measure $M_{GL,u,q}$ on
partitions, the Young Tableau Algorithm yields more: a probability
measure on standard Young tableaux. One consequence of this is a (new)
representation of prinicipally specialized Hall-Littlewood polynomials
as a sum of certain weights over standard Young tableaux.

	Let us indicate an application of this probability measure on
standard Young tableaux. Rudvalis and Shinoda $\cite{RS}$ studied the
distribution of fixed vectors for the classical groups over finite
fields. Let $G=G(n)$ be a classical group (i.e. one of $GL$,$U$,$Sp$,
or $O$) acting on an $n$ dimensional vector space $V$ over a finite
field $F_q$ (in the unitary case $F_{q^2}$) in its natural way. Let
$P_{G,n}(k,q)$ be the chance that an element of $G$ fixes a $k$
dimensional subspace and let $P_{G,\infty}(k,q)$ be the $n \rightarrow
\infty$ limit of $P_{G,n}(k,q)$. They found (in a 76 page unpublished
work) beautiful formulas for $P_{G,\infty}(k,q)$. Their formulas are
(setting $x=\frac{1}{q}$):

\begin{enumerate}
\item $P_{GL,\infty}(k,q) = \left[\prod_{r=1}^{\infty} (1-x^r)\right]
\frac{x^{k^2}}{(1-x)^2 \cdots (1-x^k)^2}$

\item $P_{U,\infty}(k,q) = \left[\prod_{r=1}^{\infty} \frac{1}{1+x^{2r-1}}\right]
\frac{x^{k^2}}{(1-x^2) \cdots (1-x^{2k})}$

\item $P_{Sp,\infty}(k,q) = \left[\prod_{r=1}^{\infty} \frac{1}{1+x^r}\right]
\frac{x^{\frac{k^2+k}{2}}}{(1-x) \cdots (1-x^k)}$

\item $P_{O,\infty}(k,q) = \left[\prod_{r=0}^{\infty}
\frac{1}{1+x^r}\right] \frac{x^{\frac{k^2-k}{2}}}{(1-x) \cdots
(1-x^k)}$.  \end{enumerate} From a probabilistic perspective, it is
very natural to try to interpret the factorizations in these formulas
as certain random variables being independent (the paper \cite{RS}
gives no insight as to why these formulas have a product form). The
Young tableau algorithm leads to such an understanding for the finite
general linear and unitary groups; see \cite{F3} for details.

{\bf Remarks}

\begin{enumerate}

\item A skew diagram is the set theoretic difference between paritions
$\mu,\lambda$ with $\mu \subseteq \lambda$ and a horizontal strip is a
skew diagram with at most one square in each column. There is another
algorithm for growing random partitions distributed according to
$M_{GL,u,q}$ in which one tosses coins and adds horizontal strips (as
opposed to a box at a time). Details are in \cite{F3}.

\item (Joint with Mark Huber) We indicate how to make the Young
Tableau Algorithm run on a computer, so as to terminate in finite time
(clearly one can't flip infinitely many coins). Let $a_N$ be the
number of times that coin $N$ comes up heads; the idea is to first
determine the random vector $(a_1,a_2,\cdots)$ and then grow the
partitions as in Step 2b of the Young Tableau Algorithm. So let us
explain how to determine $(a_1,a_2,\cdots)$. For $N \geq 1$ let
$t^{(N)}$ be the probability that all tosses of all coins numbered $N$
or greater are tails. For $N \geq 1$ and $j \geq 0$ let $t^{(N)}_j$ be
the probability that some toss of a coin numbered $N$ or greater is a
head and that coin $N$ comes up heads $j$ times. It is simple to write
down expressions for $t^{(N)},t^{(N)}_0,t^{(N)}_1,\cdots$ and clearly
$t^{(N)}+\sum_{j \geq 0} t^{(N)}_j=1$.

	The basic operation a computer can perform is to produce a
random variable $U$ distributed uniformly in the interval $[0,1]$. By
dividing $[0,1]$ into intervals of length
$t^{(1)},t^{(1)}_0,t^{(1)}_1,\cdots$ and seeing where $U$ is located,
one arrives at the value of $a_1$. Furthermore, if $U$ landed in the
interval of length $t^{(1)}$ then all coins come up tails and the
algorithm is over. Otherwise, move on to coin 2, dividing $[0,1]$ into
intervals of length $t^{(2)},t^{(2)}_0,t^{(2)}_1,\cdots$ and so on.

	For $0<u<1$ and $q$ the size of a finite field, this algorithm
terminates quickly. The probability of the algorithm stopping after
the generation of the first uniform in $[0,1]$ is
$\prod_{i=1}^{\infty} (1-u/q^i) \geq \prod_{i=1}^{\infty} (1-1/q^i) >
(1-1/q)^2 \geq 1/4$ where the second inequality is Corollary 3.6 of
\cite{NP2}. Should it be necessary to generate future uniforms, the
same argument shows that the algorithm stops after each one with
probalility at least $1/2$.

\end{enumerate}	

\subsection{Sampling for a Given Size: Unipotent Elements}
\label{samplingsize}

	An element of $GL(n,q)$ is called unipotent if all of its
eigenvalues are $1$; a theorem of Steinberg asserts that the number of
unipotent elements in $GL(n,q)$ is $q^{n(n-1)}$ (this is the square of
the order of a $q$-Sylow subgroup if $q$ is prime). Unipotent elements
are interesting because any element $\alpha$ in $GL(n,q)$ can be written
uniquely as the product $\alpha_{s} \alpha_{u}$ where $\alpha_s$ is
semisimple and $\alpha_u$ is unipotent.

	Thus it is natural to study the random partition
$\lambda_{z-1}$ for unipotent elements in $GL(n,q)$. This is the same
as conditioning the measure $M_{GL,u,q}$ to live on partitions of size
$n$. This subsection explains how to modify the sampling method of
Subsection \ref{Symmetric} to sample from this conditioned version of
$M_{GL,u,q}$ and also from a $q$-analog of Plancharel measure (related
to the longest increasing subsequence problem). These results are
joint with Mark Huber.

\begin{center}
Algorithm for Sampling from $M_{GL,u,q}$ given that $|\lambda|=n$
\end{center}

\begin{description}

\item [Step 0] Start with $N=1$ and $\lambda$ the empty partition.

\item [Step 1] If $n=0$ then stop. Otherwise set $h=1-\frac{1}{q^n}$.

\item [Step 2] Flip a coin with probability of heads $h$.

\item [Step 2a] If the toss of Step 2 came up tails, increase the value
of $N$ by $1$ and go to Step 2.

\item [Step 2b] If the toss of Step 2 comes up heads, decrease the value
of $n$ by $1$, increase $\lambda$ according to the rule of Step 2b of
the Young Tableau Algorithm (which depends on $N$), and then go to Step 1.

\end{description}

	Theorem \ref{unipotentalg} will show that the above algorithm
samples from $M_{GL,u,q}$ conditioned to live on partitions of size
$n$. It is perhaps surprising that unlike the Young Tableau Algorithm,
the probability of a coin coming up heads is independent of the coin
number; it depends only on the number of future boxes needed to get a
partition of size $n$.

\begin{lemma} \label{useful1} Let $N_i$ be the number of times that
coin $i$ comes up heads in the Young Tableau Algorithm with $u=1$ and
let $\vec{N_i}$ be the infinite vector with $i$th component $N_i$.

\begin{enumerate}

\item The probability that $\vec{N_i}=\vec{n_i}$ is
$\frac{\prod_{r=1}^{\infty} (1-\frac{1}{q^i})}{q^{\sum_{i} in_i}}$.

\item \[ \sum_{\vec{n_i}: \sum n_i=a} \frac{1}{q^{\sum_{i} in_i}} =
\frac{1}{q^a (\frac{1}{q})_a}.\]
\end{enumerate}
\end{lemma}

\begin{proof} The first assertion is clear. The second assertion is
well known in the theory of partitions, but we argue
probabilistically. Multiply both sides by $\prod_{r=1}^{\infty}
(1-\frac{1}{q^i})$. Then note from the first assertion that the left
hand side is the $M_{GL,1,q}$ chance of having a partition of size
$a$. Now use the second equation in the proof of Lemma
\ref{probabmeasure} in Subsection \ref{measures}. \end{proof}

	For Theorem \ref{unipotentalg} the notation Prob. is shorthand
for the probability of an event.

\begin{theorem} \label{unipotentalg} The algorithm for sampling from
$M_{GL,u,q}$ conditioned to live on partitions on size $n$ is
valid. \end{theorem}

\begin{proof} From the formula for $M_{GL,u,q}$, the conditioned
measure for $M_{GL,u,q}$ is the same as for $M_{GL,1,q}$. Now let
$n_i$ be the number of times that coin $i$ comes up heads in the Young
Tableau Algorithm. Letting $|$ denote conditioning, it suffices to
show that

\[ Prob. (n_i \geq 1 | \sum_{j \geq i} n_j=s) = 1-\frac{1}{q^s}.\] In
fact (for reasons to be explained later) we compute a bit more, namely
the conditional probability that $n_i=a$ given that $\sum_{j \geq i}
n_j = s$. By definition this conditional probability is the ratio

\[ \frac{Prob. (n_i=a,\sum_{j \geq i} n_j = s)}{Prob. (\sum_{j \geq i}
n_j=s)}.\] The numerator and denominator are computed using Lemma
\ref{useful2} as follows:

\begin{eqnarray*}
Prob. (n_i=a,\sum_{j \geq i} n_j = s) & = & \sum_{a_{i+1}+\cdots=s-a}
\frac{\prod_{r=i}^{\infty} (1-1/q^r)}{q^{ia} q^{\sum_{j \geq i+1} ja_j}}\\
& = &  \sum_{a_{i+1}+\cdots=s-a} \frac{\prod_{r=i}^{\infty} (1-1/q^r)}{q^{is} q^{\sum_{j \geq i+1} (j-i)a_j}}\\
& = &  \sum_{a_{1}+\cdots=s-a} \frac{\prod_{r=i}^{\infty} (1-1/q^r)}{q^{is} q^{\sum_{j \geq 1} a_j}}\\
& = & \frac{\prod_{r=i}^{\infty} (1-1/q^r)}{q^{is} q^{s-a} (\frac{1}{q})_{s-a}}.
\end{eqnarray*}

\begin{eqnarray*}
Prob. (\sum_{j \geq i} n_j=s) & = & \sum_{a_{i}+\cdots=s}
\frac{\prod_{r=i}^{\infty} (1-1/q^r)}{q^{\sum_{j \geq i} ja_j}}\\
& = & \sum_{a_{i}+\cdots=s}
\frac{\prod_{r=i}^{\infty} (1-1/q^r)}{q^{(i-1)s+\sum_{j \geq i} (j-(i-1))a_j}}\\
& = & \sum_{a_{1}+\cdots=s}
\frac{\prod_{r=i}^{\infty} (1-1/q^r)}{q^{(i-1)s+\sum_{j \geq 1} ja_j}}\\
& = & \frac{\prod_{r=i}^{\infty} (1-1/q^r)}{q^{is}(\frac{1}{q})_s}.
\end{eqnarray*} Thus $Prob. (n_i = 0 | \sum_{j \geq i} n_j=s) = \frac{1}{q^s}$
and the result follows.
\end{proof}

	As mentioned in Subsection \ref{Symmetric} there is a natural
measure $M_{Pl,q}$ on the set of all partitions of all integers which
when conditioned to live on partitions of a given size gives a
$q$-analog of Plancherel measure, which is related to longest
increasing subsequence in non-uniform random permutations
\cite{F3}. In what follows $J_a(q)$ is the polynomial discussed on
pages 52-54 of \cite{F1}, $h(s)$ denotes the hook-length of a dot in
$\lambda$ \cite{Mac} and $[n]=\frac{q^n-1}{q-1}$ is the $q$-analog of
the number $n$. Recall that a skew diagram is the set theoretic
difference between paritions $\mu,\lambda$ with $\mu \subseteq
\lambda$ and that a horizontal strip is a skew diagram with at most
one square in each column.

\begin{center}
Algorithm for Sampling from $M_{Pl,q}$ for $q > 1$ given that $|\lambda|=n$
\end{center}

\begin{description}

\item [Step 0] Start with $\lambda$ the empty partition.

\item [Step 1] If $n=0$ then stop. Otherwise choose $a$ with $0 \leq a
\leq n$ with probability

\[ \frac{q^{n^2} (1-\frac{1}{q^{n-a+1}})^2 \cdots
(1-\frac{1}{q^n})^2}{q^{(n-a)^2+n} (\frac{1}{q})_a}
\frac{J_{n-a}(q)}{J_n(q)}.\] Then increase $\lambda$ to $\Lambda$ with
probability 

\[ (1-\frac{1}{q}) \cdots (1-\frac{1}{q^a}) \frac{q^{n(\lambda)}
\prod_{s \in \lambda} (1-\frac{1}{q^{h(s)}})}{q^{n(\Lambda)} \prod_{s
\in \Lambda} (1-\frac{1}{q^{h(s)}})} \] if $\Lambda-\lambda$ is a
horizontal strip of size $a$ and with probability $0$
otherwise. Finally replace $n$ by $n-a$ and repeat Step 1.

\end{description}

	Using Lemma \ref{useful2}, Theorem \ref{unipotentalg2} proves
that the algorithm for sampling from $M_{Pl,q}$ conditioned to live on
$|\lambda|=n$ works. We omit the details, which (given the background
material in \cite{F1}) are analogous to the case of $M_{GL,u,q}$.

\begin{lemma} \label{useful2} Let $N_i$ be the number of times that
coin $i$ comes up heads in the algorithm from \cite{F3} for sampling
from the measure $M_{Pl,q}$ and let $\vec{N_i}$ be the infinite vector
with $i$th component $N_i$.
\begin{enumerate}

\item The probability that $\vec{N_i}=\vec{n_i}$ is
$\frac{\prod_{r=1}^{\infty} \prod_{j=r}^{\infty}
(1-\frac{1}{q^j})}{q^{\sum_{i} in_i} \prod_i (\frac{1}{q})_{n_i}}$.

\item \[ \sum_{\vec{n_i}: \sum_{n_i}=a} \frac{1}{q^{\sum_{i} in_i}
\prod_i (\frac{1}{q})_{n_i}} = \frac{J_a(q)}{q^{a^2} (1-\frac{1}{q})^2
\cdots (1-\frac{1}{q^a})^2}.\]
\end{enumerate}
\end{lemma}

\begin{theorem} \label{unipotentalg2} The algorithm given for sampling
from $M_{Pl,q}$ with $q > 1$ conditioned to live on $|\lambda|=n$ is
valid. \end{theorem}

\subsection{Markov Chain Approach} \label{MarkovGL}

	The main result in this subsection is a third method for
understanding the measure $M_{GL,u,q}$ probabilistically
(\cite{F7}). The idea is to build up the random partition a column at
a time; if the current column has size $a$, then the next column will
have size $b$ (with $b \leq a$) with probability $K(a,b)$. The
surprise is that this transition rule turns out to be independent of
the columns, yielding a Markov on the natural numbers. This Markov
chain is diagonalizable with eigenvalues $1,\frac{u}{q},
\frac{u^2}{q^4}, \cdots$. It will be used to give a probabilistic
proof of the Rogers-Ramanujan identities in Subsection \ref{RR}.

	It is convenient to set $\lambda_0'$ (the height of an
imaginary zeroth column) equal to $\infty$. For the entirety of this
subsection, let $P(a)$ be the $M_{GL,u,q}$ probability that
$\lambda_1'=a$. Theorem \ref{markovgl}, which makes the connection
with Markov chains, is proved in a completely elementary way. The
argument reproves that $M_{GL,u,q}$ is a probability measure (Lemma
\ref{probabmeasure} of Subsection \ref{measures}), shows that the
asserted Markov transition probabilities add to one, and gives a
formula for $P(a)$.

\begin{theorem} \label{markovgl} Starting with $\lambda_0'=\infty$,
define in succession $\lambda_1',\lambda_2',\cdots$ according to the
rule that if $\lambda_i'=a$, then $\lambda_{i+1}'=b$ with probability

\[ K(a,b) = \frac{u^b (\frac{1}{q})_a (\frac{u}{q})_a}{q^{b^2}
(\frac{1}{q})_{a-b} (\frac{1}{q})_b (\frac{u}{q})_b}.\] Then the
resulting partition is distributed according to
$M_{GL,u,q}$. \end{theorem}

\begin{proof} Suppose we know that $M_{GL,u,q}$ is a probability
measure and that \[ P(a) = \frac{u^a (\frac{u}{q})_{\infty}}{q^{a^2}
(\frac{1}{q})_a (\frac{u}{q})_a}. \] Then the $M_{GL,u,q}$ probability
of choosing a partition with $\lambda_i'=r_i'$ for all $i$ is

\[ Prob.(\lambda_0'=\infty) \frac{Prob.(\lambda_0'=\infty,
\lambda_1'=r_1)} {Prob.(\lambda_0'=\infty)} \prod_{i=1}^{\infty}
\frac{Prob.(\lambda_0'=\infty,
\lambda_1'=r_1,\cdots,\lambda_{i+1}'=r_{i+1})}
{Prob.(\lambda_0'=\infty,\lambda_1'=
r_1,\cdots,\lambda_{i}'=r_{i})}.\] Thus it is enough to prove the
(surprising) assertion that

\[ \frac{Prob.(\lambda_0'=\infty,
\lambda_1'=r_1,\cdots,\lambda_{i-1}'=r_{i-1}, \lambda_i'=a,
\lambda_{i+1}'=b)} {Prob.(\lambda_0'=\infty,\lambda_1'=
r_1,\cdots,\lambda_{i-1}'=r_{i-1}, \lambda_{i}'=a)} = \frac{u^b
(\frac{1}{q})_a (\frac{u}{q})_a}{q^{b^2} (\frac{1}{q})_{a-b}
(\frac{1}{q})_b (\frac{u}{q})_b},\] for all
$i,a,b,r_1,\cdots,r_{i-1}$. One calculates that

\[ \sum_{\lambda: \lambda_1'=r_1,\cdots,\lambda_{i-1}'=r_{i-1} \atop
\lambda_i'=a} M_{GL,u,q}(\lambda) = \frac{u^{r_1+\cdots+r_{i-1}}}
{q^{r_1^2+\cdots+r_{i-1}^2} (\frac{1}{q})_{r_1-r_2} \cdots
(\frac{1}{q})_{r_{i-2}-r_{i-1}} (\frac{1}{q})_{r_{i-1}-a}} P(a).\]
Similarly, observe that \[ \sum_{\lambda:
\lambda_1'=r_1,\cdots,\lambda_{i-1}'=r_{i-1} \atop
\lambda_i'=a,\lambda_{i+1}'=b} M_{GL,u,q}(\lambda) =
\frac{u^{r_1+\cdots+r_{i-1}+a}} {q^{r_1^2+\cdots+r_{i-1}^2+a^2}
(\frac{1}{q})_{r_1-r_2} \cdots (\frac{1}{q})_{r_{i-2}-r_{i-1}}
(\frac{1}{q})_{r_{i-1}-a} (\frac{1}{q})_{a-b}} P(b).\] Thus the ratio
of these two expressions is \[ \frac{u^b (\frac{1}{q})_a
(\frac{u}{q})_a}{q^{b^2} (\frac{1}{q})_{a-b} (\frac{1}{q})_b
(\frac{u}{q})_b}, \] as desired. Note that the transition
probabilities must sum to 1 because

\[ \sum_{b \leq a} \frac{ \sum_{\lambda:
\lambda_1'=r_1,\cdots,\lambda_{i-1}'=r_{i-1} \atop
\lambda_i'=a,\lambda_{i+1}'=b} M_{GL,u,q}(\lambda)}{\sum_{\lambda:
\lambda_1'=r_1,\cdots,\lambda_{i-1}'=r_{i-1} \atop \lambda_i'=a}
M_{GL,u,q}(\lambda)} = 1\] for any measure $M_{GL,u,q}$ on
partitions.

	Thus to complete the proof, it must be shown that $M_{GL,u,q}$
is a probability measure and that \[ P(a) = \frac{u^a
(\frac{u}{q})_{\infty}}{q^{a^2} (\frac{1}{q})_a (\frac{u}{q})_a}. \]
Since \[ \frac{ \sum_{\lambda:
\lambda_1'=r_1,\cdots,\lambda_{i-1}'=r_{i-1} \atop
\lambda_i'=a,\lambda_{i+1}'=b} M_{GL,u,q}(\lambda)}{\sum_{\lambda:
\lambda_1'=r_1,\cdots,\lambda_{i-1}'=r_{i-1} \atop \lambda_i'=a}
M_{GL,u,q}(\lambda)} = \frac{P(b)u^a}{P(a) q^{a^2}
(\frac{1}{q})_{a-b}} \] it follows that \[ \sum_{b \leq a}
\frac{P(b)u^a}{P(a) q^{a^2} (\frac{1}{q})_{a-b}}=1.\] From this
recursion and the fact that $P(0)=(\frac{u}{q})_{\infty}$, one solves
for $P(a)$ inductively, finding that \[ P(a) = \frac{u^a
(\frac{u}{q})_{\infty}} {q^{a^2} (\frac{1}{q})_a (\frac{u}{q})_a} .\]
Cauchy's identity (page 20 of \cite{A1}) gives that $\sum_a
P(a)=1$, so that $M_{GL,u,q}$ is a probability measure. \end{proof}

	Theorem \ref{diagonalize} diagonalizes the transition matrix
$K$, finding a basis of eigenvectors, which is fundamental for
understanding the Markov chain (part 3 is stated as a Lemma in
\cite{A2}). Since the matrix $K$ is upper triangular with distinct
eigenvalues, this is straightforward.

\begin{theorem} \label{diagonalize} \begin{enumerate}

\item Let $C$ be the diagonal matrix with $(i,i)$ entry
$(\frac{1}{q})_i (\frac{u}{q})_i$. Let $M$ be the matrix $\left(
\frac{u^j}{q^{j^2} (\frac{1}{q})_{i-j}} \right)$. Then $K=CMC^{-1}$,
which reduces the problem of diagonalizing $K$ to that of
diagonalizing $M$.

\item Let $A$ be the matrix $\left( \frac{1}{(\frac{1}{q})_{i-j}
(\frac{u}{q})_{i+j}} \right)$. Then the columns of $A$ are
eigenvectors of $M$ for right multiplication, the $j$th column having
eigenvalue $\frac{u^j}{q^{j^2}}$.

\item The inverse matrix $A^{-1}$ is $\left( \frac{(1-u/q^{2i})
(-1)^{i-j} (\frac{u}{q})_{i+j-1}}{q^{i-j \choose 2}
(\frac{1}{q})_{i-j}} \right)$.

\end{enumerate}
\end{theorem}

	Corollary \ref{bigcor} (immediate from Theorem
\ref{diagonalize}) will be useful for the proof of the
Rogers-Ramanujan identities in Section \ref{RR}. In the case $L
\rightarrow \infty$ and $j=0$, it is the so called Rogers-Selberg
identity.

\begin{cor} \label{bigcor} Let $E$ be the diagonal matrix with $(i,i)$
entry $\frac{u^i}{q^{i^2}}$. Then $K^r=C A E^r A^{-1} C^{-1}$. More
explicitly,

\[ K^r(L,j) = \frac{(\frac{1}{q})_L (\frac{u}{q})_L}{(\frac{1}{q})_j
(\frac{u}{q})_j} \sum_{n=0}^{\infty} \frac{u^{rn} (1-u/q^{2n})
(-1)^{n-j} (\frac{u}{q})_{n+j-1}}{q^{rn^2} (\frac{1}{q})_{L-n}
(\frac{u}{q})_{L+n} q^{n-j \choose 2} (\frac{1}{q})_{n-j}}.\] \end{cor}

\begin{proof} This is immediate from Theorem
\ref{diagonalize}. \end{proof}

{\bf Remarks}

\begin{enumerate}

\item One of our motivations for seeking a Markov chain description of
$M_{GL,u,q}$ is work of Fristedt \cite{Fr}, who had a Markov chain
approach for the measure $P_q$ on the set of all partitions of all
natural numbers defined by $P_q(\lambda) = \prod_{i=1}^{\infty}
(1-q^i) q^{|\lambda|}$ where $q<1$. Fristedt's interest was in
studying what a uniformly chosen partition of an integer looks like,
and conditioning $P_q$ to live on partitions of size $n$ gives a
uniform partition. The measure $P_q$ is related to to vertex operators
\cite{O1} and to the enumeration of ramified coverings of the torus
\cite{D}. In this regard the papers \cite{O1} and \cite{BlO} prove that
the $k$ point correlation function

\[ F(t_1,\cdots,t_k)= \sum_{\lambda} q^{|\lambda|} \prod_{k=1}^n
\sum_{i=1}^{\infty} t_k^{\lambda_i-i+\frac{1}{2}}\] is a sum of
determinants involving genus 1 theta functions and their derivatives
and give connections with quasi-modular forms. It would be marvellous
if the measure $M_{GL,u,q}$ (being related to modular forms via the
Rogers-Ramanujan identities) is also related to enumerative questions
in algebraic geometry.

\item As mentioned in the introduction, the Markov chain approach
gives a unified description of conjugacy classes of the finite
classical groups. For the symplectic and orthogonal groups it is
necessary to use two Markov chains $K_1$ and $K_2$. For the symplectic
case, steps with column number $i$ odd use $K_1$ and steps with column
number $i$ even use $K_2$. For the orthogonal case, steps with column
number $i$ odd use $K_2$ and steps with column number $i$ even use
$K_1$. The Markov chains $K_1,K_2$ are the same for both cases!
Details are in \cite{F6}. The Markov chain approach is also related to
quivers \cite{F7}.

\end{enumerate}

\subsection{Rogers-Ramanujan Identities} \label{RR}

	The Rogers-Ramanujan identities \cite{Ro}

\[ 1+\sum_{n=1}^{\infty} \frac{q^{n^2}}{(1-q)(1-q^2) \cdots (1-q^n)} =
\prod_{n=1}^{\infty} \frac{1}{(1-q^{5n-1})(1-q^{5n-4})} \]

\[ 1+\sum_{n=1}^{\infty} \frac{q^{n(n+1)}}{(1-q)(1-q^2) \cdots
(1-q^n)} = \prod_{n=1}^{\infty} \frac{1}{(1-q^{5n-2})(1-q^{5n-3})} \]
are among the most interesting partition identities in number theory
and combinatorics, with connections to Lie theory and statistical
mechanics (see the discussions in \cite{A2} and \cite{F7} for many
references). One ongoing challenge in the subject (posed by Hardy) has
been to find a proof of the Rogers-Ramanujan identities which is both
motivated and simple. The purpose of this subsection is to describe
such a proof (\cite{F7}), which is also the first probabilistic proof
of the Rogers-Ramanujan identities.

	To illustrate the idea we give the proof of the following
generalization of the first Rogers-Ramanujan identity (called the
Andrews-Gordon identity \cite{A3},\cite{Go}):

\[ \sum_{n_1,\cdots,n_{k-1} \geq 0} \frac{1}{q^{N_1^2 + \cdots +
N_{k-1}^2}(1/q)_{n_1}\cdots(1/q)_{n_{k-1}}} = \prod_{r=1 \atop r \neq
0, \pm k (mod \ 2k+1)}^{\infty} \frac{1}{1-(1/q)^r} \] where
$N_i=n_i+\cdots+n_{k-1}$.

	The idea is simple. We study the distribution of the length of
the first row of a random partition distributed as $M_{GL,1,q}$. From
the definition of $M_{GL,1,q}$ the probability that the first row has
length less than $k$ is equal to

\[ \prod_{r=1}^{\infty} (1-\frac{1}{q^r}) \sum_{\lambda: \lambda_k'=0}
\frac{1}{q^{(\lambda_1')^2 +\cdots+ (\lambda_{k-1}')^2}
(1/q)_{\lambda_1'-\lambda_2'} \cdots
(1/q)_{\lambda_{k-1}'-\lambda_{k}'}} .\] Letting $n_i$ denote
$\lambda_{i}'-\lambda_{i+1}'$ and $N_i$ denote $\lambda_i'$, this becomes

\[ \prod_{r=1}^{\infty} (1-\frac{1}{q^r}) \sum_{n_1,\cdots,n_{k-1}
\geq 0} \frac{1}{q^{N_1^2 + \cdots +
N_{k-1}^2}(1/q)_{n_1}\cdots(1/q)_{n_{k-1}}} \] which is a essentially
the left hand side of the Andrews-Gordon identity. On the other hand
the probability that the first row has length less than $k$ is equal
to the probability that the Markov chain of Section \ref{MarkovGL} is
absorbed at $0$ at time $k$. Since we diagonalized the matrix
associated to this Markov chain, it is straightforward to compute this
probability. To get it into product form it is necessary to apply
Jacobi's triple product identity which has a simple combinatorial
proof \cite{A1}. Further details are in \cite{F7}.

	Next we argue that this proof is motivated. Certainly the
measure $M_{GL,1,q}$ is a natural object to study, given that it is
the $n \rightarrow \infty$ limit law of $\lambda_{z-1}$ for a random
element of $GL(n,q)$. It was natural to try to build up the random
partitions $\lambda$ column by column as in Section
\ref{MarkovGL}. Observing that the resulting Markov chain is absorbing
at $0$ with probability one, the time to absorption (equivalent to the
distribution of the length of the first row) is the most natural
quantity one could examine. The final step is applying Jacobi's triple
product identity, and thus going from a ``sum = sum'' identity to a
``sum = product'' identity. As mentioned above Jacobi's triple product
identity is easy to verify, but one still wants a motivation for
trying to write the left hand side of the Andrews-Gordon identity in
product form. One motivation is Baxter's work on statistical mechanics
(surveyed in \cite{A2},\cite{Bax1},\cite{Bax2}) in which he really
needed ``sum = product'' identities and was led to conjecture analogs
of Rogers-Ramanujan type identities. Although a proof of the
Rogers-Ramanujan identities doesn't emerge from his work, it is
clearly one of the truly great accomplishments in mathematics and his
book \cite{Bax1} has been very influential. A second motivation is our
work on the $n \rightarrow \infty$ asymptotic probability that an
element of $GL(n,q)$ is semisimple. The argument, recorded in
\cite{F1} or the more readily available \cite{F4} needed a ``sum =
product'' identity. The corresponding computation in \cite{F9} for the
finite affine groups needed both Rogers-Ramanujan identities.

	Andrews' paper \cite{A4} notes that many proofs of the
Rogers-Ramanujan identities make use of the following result called
Bailey's Lemma, alluded to in \cite{Ba} and stated explicity in
\cite{A3}. A pair of sequences $\{\alpha_L\}$ and $\{\beta_L\}$ are
called a Bailey pair if \[ \beta_L = \sum_{r=0}^L
\frac{\alpha_r}{(1/q)_{L-r} (u/q)_{L+r}}.\] Bailey's Lemma states that
if $\alpha_L'=\frac{u^L}{q^{L^2}} \alpha_L$ and $\beta_L'=\sum_{r=0}^L
\frac{u^r}{q^{r^2}(1/q)_{L-r}} \beta_r$, then $\{\alpha_L'\}$ and
$\{\beta_L'\}$ are a Bailey pair. From the viewpoint of Markov chains,
this case of Bailey's Lemma is clear. To explain, let $A,D,M$ be as in
Theorem \ref{diagonalize} (recall that $M=ADA^{-1}$). Viewing
$\alpha=\vec{\alpha_L}$ and $\beta=\vec{\beta_L}$ as column vectors,
the notion of a Bailey pair means that $\beta=A \alpha$. This case of
Bailey's Lemma follows because \[ \beta' = M \beta = ADA^{-1} \beta =
AD \alpha = A \alpha'.\] As Andrews explains in \cite{A2}, the power
of Bailey's lemma lies in its ability to be iterated and gives a short
proof of the Rogers-Selberg identity (Corollary \ref{bigcor} in
Section \ref{MarkovGL}). From the remarks in this paragraph it is
clear that iterating Bailey's lemma corresponds to taking several
according to the Markov chain $K$. This demystifies the Bailey's Lemma
proofs of the Rogers-Ramanujan identities, which strike this author as
unmotivated. The fact that the Markov chain approach has analogs for
other finite classical groups and for quivers is further evidence of
its naturality.

\section{Upper Triangular Matrices} \label{triangular}

	This section surveys probabilistic aspects of conjugacy
classes in the group $T(n,q)$ of upper triangular matrices over finite
fields with $1$'s along the main diagonal. At present little is known
about conjugacy in $T(n,q)$. The papers \cite{VAr}, \cite{VArV} study
the number of conjugacy classes. Kirillov \cite{Kir} calls for an
extension of his method of coadjoint orbits for groups over real,
complex, or $p$-adic fields to the group $T(n,q)$ and gives
premilinary connections with statistical physics; the paper \cite{IK}
gives a counterexample to one of his conjectures. As we do not see how
to further develop those results or improve on their exposition, we
instead focus on a simpler problem: the probabilistic study of Jordan
form of elements of $T(n,q)$.

	Subsection \ref{growth} describes a probabilistic growth
algorithm for the Jordan form of upper triangular matrices over a
finite field. This is linked with symmetric function theory and
potential theory on Bratteli diagrams in Subsection \ref{Bratteli}.

\subsection{Growth Algorithm for Jordan Form} \label{growth}

	Theorem \ref{growthtriangular} gives a probabilistic growth
algorithm for the Jordan form of random elements of $T(n,q)$. Its
proof uses elementary reasoning from linear algebra.

\begin{theorem} \label{growthtriangular} (\cite{Kir},\cite{B}) The
Jordan form of a uniformly chosen element of $T(n,q)$ can be sampled
from by stopping the following procedure after $n$ steps:

	Starting with the empty partition, at each step transition
from a partition $\lambda$ to a partition $\Lambda$ by adding a box to
column $i$ chosen according to the rules

\begin{itemize} \item $ i=1$ with probability
$\frac{1}{q^{\lambda_1'}}$ \item $i=j>1$ with probability
$\frac{1}{q^{\lambda_j'}}-\frac{1}{q^{\lambda_{j-1}'}}$ 
\end{itemize}
\end{theorem}

	Theorem \ref{growthtriangular} leads to the following central
limit theorem about the asymptotic Jordan form of an element of
$T(n,q)$.

\begin{theorem} (\cite{B}) Let $\lambda$ be the partition
corresponding to the Jordan form of a random element of $T(n,q)$. Let
$Prob^n$ denote probability under the uniform measure on $T(n,q)$ and
let $p_i=\frac{1}{q^{i-1}}-\frac{1}{q^i}$. Then

\[ lim_{n \rightarrow \infty} Prob^{n} (\frac{\lambda_i-p_i n}{\sqrt{n}} \leq x_i, i=1, \cdots, k ) = (2
\pi)^{-\frac{k}{2}} \int_{-\infty}^{x_1} \cdots \int_{\infty}^{x_k}
e^{-\frac{1}{2} <Qt,t>} dt \] for any $(x_1,\cdots,x_k) \in R^k$,
where the covariance matrix equals

\[ Q = diag(p_1,\cdots,p_k) - (p_ip_j)_{i,j=1}^k.\] \end{theorem}

\subsection{Symmetric Functions and Potential Theory} \label{Bratteli}

	Given the usefulness of symmetric functions in the
probabilistic study of the measure $M_{GL,u,q}$, it is natural to seek
an analogous understanding of Theorem \ref{growthtriangular}. That is
the topic of the present subsection. The ideas here are from the
report \cite{F8}.

	The first step is to link the probability that an element of
$T(n,q)$ has Jordan form of type $\Lambda$ with symmetric function
theory. For the rest of this section, $P_{\Lambda}(q,t)$ denotes a
Macdonald polynomial, $K_{\mu \Lambda}(q,t)$ denotes a Kostka-Foulkes
polynomial, and $f^{\mu}$ is the dimension of the irreducible
representation of $S_n$ corresponding to the partition $\mu$ (see
\cite{Mac} for background). Note that when $q=0$ the Macdonald
polynomial is our friend, a Hall-Littlewood polynomial.

\begin{theorem} \label{formup} (\cite{F5}) The probability that a
random element of $T(n,q)$ has Jordan form of type $\Lambda$ is

\[ P_{\Lambda}(1-\frac{1}{q},\frac{1}{q}
-\frac{1}{q^2},\cdots;0,\frac{1}{q}) \sum_{\mu \vdash n} f^{\mu}
K_{\mu \Lambda}(0,q).\] \end{theorem}

	Next we give some background on potential theory on Bratteli
diagrams. This is a beautiful subject, with connections to probability
and representation theory. We recommend \cite{Ke} for background on
potential theory with many examples and \cite{BO} for a survey of
recent developments. The basic set-up is as follows. One starts with a
Bratteli diagram; that is an oriented graded graph $\Gamma= \cup_{n
\geq 0} \Gamma_n$ such that

\begin{enumerate}
\item $\Gamma_0$ is a single vertex $\emptyset$.
\item If the starting vertex of an edge is in $\Gamma_i$, then its end vertex is in $\Gamma_{i+1}$.
\item Every vertex has at least one outgoing edge.
\item All $\Gamma_i$ are finite.
\end{enumerate}

For two vertices $\lambda, \Lambda \in \Gamma$, one writes $\lambda
\nearrow \Lambda$ if there is an edge from $\lambda$ to
$\Lambda$. Part of the underlying data is a multiplicity function
$\kappa(\lambda,\Lambda)$. Letting the weight of a path in $\Gamma$ be
the product of the multiplicities of its edges, one defines the
dimension $dim(\Lambda)$ of a vertex $\Lambda$ to be the sum of the
weights over all maximal length paths from $\emptyset$ to $\Lambda$
(this definition clearly extend to intervals). Given a Bratteli
diagram with a multiplicity function, one calls a function $\phi$ {\it
harmonic} if $\phi(0)=1$, $\phi(\lambda) \geq 0$ for all $\lambda \in
\Gamma$, and \[ \phi(\lambda) = \sum_{\Lambda: \lambda \nearrow
\Lambda} \kappa(\lambda,\Lambda) \phi(\Lambda).\] An equivalent
concept is that of coherent probability distributions. Namely a set
$\{M_n\}$ of probability distributions $M_n$ on $\Gamma_n$ is called
{\it coherent} if \[ M_{n-1}(\lambda) = \sum_{\Lambda: \lambda
\nearrow \Lambda} \frac{dim(\lambda)
\kappa(\lambda,\Lambda)}{dim(\Lambda)} M_{n}(\Lambda).\] The formula
allowing one to move between the definitions is $\phi(\lambda) =
\frac{M_n(\lambda)}{dim(\lambda)}$.

	One reason the set-up is interesting from the viewpoint of
probability theory is the fact that every harmonic function can be
written as a Poisson integral over the set of extreme harmonic
functions (which is often the Martin boundary). For the Pascal lattice
(vertices of $\Gamma_n$ are pairs $(k,n)$ with $k=0,1,\cdots,n$ and
$(k,n)$ is connected to $(k,n+1)$ and $(k+1,n+1)$), this fact is the
simplest instance of de Finetti's theorem. When the multiplicity
function $\kappa$ is integer valued, one can define a sequence of
algebras $A_n$ associated to the Bratteli diagram, and harmonic
functions correspond to certain characters of the inductive limit of
the algebras $A_n$.

	Next we define a branching for which the probability that an
element of $T(n,q)$ has Jordan type $\Lambda$ is a harmonic
function. First some notation is needed. For $\lambda \nearrow
\Lambda$, let $R_{\Lambda / \lambda}$ (resp. $C_{\Lambda / \lambda}$)
be the boxes of $\lambda$ in the same row (resp. column) as the boxes
removed from $\lambda$ to get $\Lambda$. This notation differs from
that in \cite{Mac}. Let $a_{\lambda}(s)$, $l_{\lambda}(s)$ be the
number of dots in $\lambda$ strictly to the east and south of $s$, and
let $h_{\lambda}(s)=a_{\lambda}(s)+l_{\lambda}(s)+1$.

{\bf Definition 1:} For $0 \leq q <1$ and $0<t<1$, the underlying
Bratteli diagram $\Gamma$ has as level $\Gamma_n$ all partitions
$\lambda$ of $n$. Letting $i$ be the column number of the dot removed
to go from $\lambda$ to $\Lambda$, for $\lambda \nearrow \Lambda$,
define the multiplicty function as

\[ \kappa(\lambda,\Lambda) = \frac{1}{t^{\Lambda_i'-1}} \prod_{s \in
R_{\Lambda / \lambda}} \frac{1-q^{a_{\Lambda}(s)+1}
t^{l_{\Lambda}(s)}} {1-q^{a_{\lambda}(s)+1}t^{l_{\lambda}(s)}}
\prod_{s \in C_{\Lambda / \lambda}}
\frac{1-q^{a_{\Lambda}(s)}t^{l_{\Lambda}(s)+1}}
{1-q^{a_{\lambda}(s)}t^{l_{\lambda}(s)+1}}.\]

	Equation I.10 of \cite{GH} proves that

\[ dim(\Lambda)=\frac{1}{t^{n(\Lambda)}} \sum_{\mu \vdash n} f^{\mu}
K_{\mu \Lambda} (q,t).\]

{\bf Definition 2:} For $0 \leq q <1, 0<t<1$ and $0 \leq
x_1,x_2,\cdots$ such that $\sum x_i=1$, define a family $\{M_n\}$ of
probability measures on partitions of size $n$ by

\begin{eqnarray*}
M_n(\Lambda) & = & \frac{(1-q)^{|\Lambda|} P_{\Lambda}(x;q,t)
\sum_{\mu \vdash n} f^{\mu} K_{\mu \Lambda}(q,t)} {\prod_{s \in
\Lambda} (1-q^{a_{\Lambda}(s)+1} t^{l_{\Lambda}(s)})}\\
& = & \frac{(1-q)^{|\Lambda|} P_{\Lambda}(x;q,t) t^{n(\Lambda)} dim(\Lambda)} {\prod_{s \in
\Lambda} (1-q^{a_{\Lambda}(s)+1} t^{l_{\Lambda}(s)})}
\end{eqnarray*}

	Consider the specialization that $q=0$ and $t=\frac{1}{q}$,
where this second $q$ is the size of a finite field. Further, set
$x_i=\frac{1} {q^{i-1}}-\frac{1}{q^i}$. Then Theorem \ref{formup}
implies that $M_n(\Lambda)$ is the probability that a uniformly chosen
element of $T(n,q)$ has Jordan type $\Lambda$. The multiplicities have
a simple description; letting $i$ be the column to which one adds in
order to go from $\lambda$ to $\Lambda$, it follows that
$\kappa(\lambda,\Lambda)=q^{\lambda_i'}+ q^{\lambda_i'-1}+ \cdots+
q^{\lambda_{i+1}'}$. Second, $dim(\Lambda)$ reduces to a Green's
polynomial $Q^{\Lambda}(q)=Q^{\Lambda}_{(1^n)}(q)$ as in Section 3.7
of \cite{Mac}. These polynomials are important in the representation
theory of the finite general linear groups. This specialization was
the motivation for Definition 2.

	The connection with potential theory is given by the following
result.

\begin{theorem} (\cite{F8}) The measures of Definition 2 are harmonic
with respect to the branching of Definition 1. \end{theorem}

	It is elementary and well-known that if one starts at the
empty partition and transitions from $\lambda$ to $\Lambda$ with
probability $\frac{\kappa(\lambda,\Lambda) M_n(\Lambda) dim(\lambda)}
{M_{n-1}(\lambda) dim(\Lambda)}$, one gets samples from any coherent
family of measures $\{M_n\}$. Applying this principle to the above
specialization in which $M_n(\Lambda)$ is $T(n,q)$ and using
Macdonald's principal specialization formula (page 337 of \cite{Mac})
gives the advertised proof of Theorem \ref{growthtriangular} by means
of symmetric functions and potential theory.

{\bf Remarks:}
\begin{enumerate}

\item As indicated in \cite{F8}, the example of Schur functions
($q=t<1$) is also interesting. The measure $M_n(\Lambda)$ reduces to
$s_{\Lambda} f^{\Lambda}$, where $s_{\Lambda}$ is a Schur
function. Setting $x_1=\cdots=x_n=\frac{1}{n}$ and letting $n
\rightarrow \infty$, one obtains Plancherel measure, which is
important in representation theory and random matrix theory. Letting
$x_1=\cdots=x_n$ satisfy $\sum x_i=1$ (all other $x_j=0$) gives a
natural deformation of Plancherel measure, studied for instance by
\cite{ITW}. Stanley \cite{Sta2} shows that this measure on partitions
also arises by applying the Robinson-Schensted-Knuth algorithm to a
random permutation distributed after a biased riffle shuffle (in other
words, this measure encodes information about the longest increasing
subsequences of permutations distributed as shuffles).

\item It has been pointed out to the author that the branchings
$\kappa(\lambda,\Lambda)$ of Definition 1 are related to the
branchings $\tau(\lambda,\Lambda)$ of \cite{Ke2} by the formula \[
\kappa(\lambda,\Lambda) = f(\lambda) \tau(\lambda,\Lambda)
f(\Lambda)^{-1},\] for a certain positive function $f(\lambda)$ on the
set of vertices, which implies by \cite{Ke3} that the boundaries of
these two branchings are homeomorphic and that the branchings of
Definition 1 are multiplicative. Kerov \cite{Ke2} has a conjectural
description of the boundary. It has been verified for Schur functions
\cite{T}, Kingman branching \cite{Kin}, and Jack polynomials
\cite{KOO}, but remains open for the general case of Macdonald
polynomials. In particular, it is open for Hall-Littlewood
polynomials, the case related to $T(n,q)$. It is interesting that the
$\kappa(\lambda,\Lambda)$ of Definition 1 are integers for
Hall-Littlewood polynomials, whereas the $\tau(\lambda,\Lambda)$ of
\cite{Ke2} are not.

\end{enumerate}

\section*{Acknowledgements} The author's greatest thanks go to Persi
Diaconis (his former thesis advisor) for years of friendship,
encouragement, and inspiration. He was very helpful in the preparation
of this article. We thank Peter M. Neumann and Cheryl E. Praeger for
countless conversations about conjugacy classes and computational
group theory, and Mark Huber for permission to survey some joint
unpublished results. The author received the financial support of an
NSF Postdoctoral Fellowship.

\end{document}